\documentclass[A4paper]{article}
\usepackage{amsthm,amssymb,amsmath,color,natbib,pdfsync}

\newcommand{\levy}{L\'evy}

\newcommand{\clt}{central limit theorem}

\newcommand{\garch}{{\rm GARCH}$(1,1)$}

\newcommand{\sta}{St\u aric\u a}
\newcommand{\ex}{{\rm e}\,}
\newcommand{\sas}{s$\alpha$s}

\newcommand{\asy}{asymptotic}

\newcommand{\ts}{time series}

\definecolor{darkblue}{rgb}{.1, 0.1,.8}
\definecolor{darkgreen}{rgb}{0,0.8,0.2}
\definecolor{darkred}{rgb}{.8, .1,.1}

\def\tag{\refstepcounter{equation}\leqno }

\newtheorem{lemma}{Lemma}[section]

\newtheorem{theorem}[lemma]{Theorem}

\newtheorem{proposition}[lemma]{Proposition}
\newtheorem{definition}[lemma]{Definition}
\newtheorem{corollary}[lemma]{Corollary}
\newtheorem{example}[lemma]{Example}
\newtheorem{exercise}[lemma]{Exercise}
\newtheorem{remark}[lemma]{Remark}
\newtheorem{fig}[lemma]{Figure}
\newtheorem{tab}[lemma]{Table}

\newcommand{\MC}{Markov chain}

\newcommand{\bth}{\begin{theorem}}
\newcommand{\ethe}{\end{theorem}}
\newcommand{\sv}{stochastic volatility}
\newcommand{\bre}{\begin{remark}\em }
\newcommand{\ere}{\end{remark}}

\newcommand{\ble}{\begin{lemma}}
\newcommand{\ele}{\end{lemma}}
\newcommand{\sre}{stochastic recurrence equation}
\newcommand{\pp}{point process}
\newcommand{\bde}{\begin{definition}}
\newcommand{\ede}{\end{definition}}
\newcommand{\bco}{\begin{corollary}}
\newcommand{\eco}{\end{corollary}}

\newcommand{\bpr}{\begin{proposition}}
\newcommand{\epr}{\end{proposition}}

\newcommand{\bexer}{\begin{exercise}}
\newcommand{\eexer}{\end{exercise}}

\newcommand{\bexam}{\begin{example}}
\newcommand{\eexam}{\end{example}}

\newcommand{\bfi}{\begin{fig}}
\newcommand{\efi}{\end{fig}}

\newcommand{\btab}{\begin{tab}}
\newcommand{\etab}{\end{tab}}

\newcommand{\lhs}{left-hand side}
\newcommand{\fidi}{finite-dimensional distribution}
\newcommand{\rv}{random variable}

\newcommand{\sign}{{\rm sign}}

\newcommand{\var}{{\rm var}}

\newcommand{\cov}{{\rm cov}}

\newcommand{\rhs}{right-hand side}

\newcommand{\beao}{\begin{eqnarray*}}
\newcommand{\eeao}{\end{eqnarray*}\noindent}

\newcommand{\beam}{\begin{eqnarray}}
\newcommand{\eeam}{\end{eqnarray}\noindent}

\newcommand{\beqq}{\begin{equation}}
\newcommand{\eeqq}{\end{equation}\noindent}

\newcommand{\bce}{\begin{center}}
\newcommand{\ece}{\end{center}}

\newcommand{\barr}{\begin{array}}
\newcommand{\earr}{\end{array}}

\newcommand{\stp}{\stackrel{P}{\rightarrow}}
\newcommand{\std}{\stackrel{d}{\rightarrow}}
\newcommand{\stas}{\stackrel{\rm a.s.}{\rightarrow}}

\newcommand{\stv}{\stackrel{v}{\rightarrow}}

\newcommand{\eqd}{\stackrel{d}{=}}

\newcommand{\vague}{\stackrel{\lower0.2ex\hbox{$\scriptscriptstyle
                    \it{v} $}}{\rightarrow}}
\newcommand{\weak}{\stackrel{\lower0.2ex\hbox{$\scriptscriptstyle
                    \it{w} $}}{\rightarrow}}
\newcommand{\what}{\stackrel{\lower0.2ex\hbox{$\scriptscriptstyle
                    \it{\hat{w}} $}}{\rightarrow}}

\newcommand{\bdis}{\begin{displaymath}}
\newcommand{\edis}{\end{displaymath}\noindent}

\newcommand{\R}{\mathbb{R}}

\newcommand{\nto}{n\to\infty}

\newcommand{\xto}{x\to\infty}

\newcommand{\ov}{\overline}
\newcommand{\wt}{\widetilde}

\newcommand{\vep}{\varepsilon}

\newcommand{\la}{\lambda}

\newcommand{\halmos}{\hfill \qed}
\newcommand{\regvary}{regularly varying}
\newcommand{\slvary}{slowly varying}
\newcommand{\regvar}{regular variation}

\newcommand{\bbr}{{\mathbb R}}

\newcommand{\bbz}{{\mathbb Z}}

\newcommand{\con}{convergence}

\newcommand{\evt}{extreme value theory}

\newcommand{\st}{such that}
\newcommand{\fif}{if and only if}

\newcommand{\chf}{characteristic function}
\newcommand{\fct}{function}

\newcommand{\ds}{distribution}

\newcommand{\rep}{representation}

\newcommand{\seq}{sequence}

\newcommand{\pro}{probabilit}

\newcommand{\ms}{measure}

\newcommand{\ld}{large deviation}
\newcommand{\bfx}{{\bf x}}

\def\1{\ensuremath{\mathrm{1}\hspace{-.35em} \mathrm{1}}} 

\def\E{{\mathbb E}}

\def\P{{\mathbb{P}}}
\def\R{\mathbb{R}}
\def\Z{\mathbb{Z}}

\renewcommand{\le}{\ensuremath{\leqslant}}
\renewcommand{\ge}{\ensuremath{\geqslant}}

\newcommand{\introo}[2]{{\left]{#1,\,#2\,}\right[\kern1pt}}

\newcommand{\intrfo}[2]{{\left[{#1,\,#2}\right[\kern1pt}}

\begin{document}

\title{Stable limits for sums of dependent infinite variance random variables}


\author{
Katarzyna Bartkiewicz  and Adam Jakubowski \\\small Nicolaus Copernicus University,\\\small
            Faculty of Mathematics and Computer Science \\\small
              ul. Chopina 12/18,
              87-100 Toru\'n, Poland\\\small
             \texttt{kabart@mat.umk.pl} and \texttt{adjakubo@mat.umk.pl}\and
Thomas Mikosch\footnote{Thomas Mikosch's research is partly supported
by the Danish Research Council (FNU) Grants 272-06-0442 and 09-072331. The
research of Thomas Mikosch and Olivier Wintenberger is partly supported by a
Scientific Collaboration Grant of the French Embassy in Denmark.}  \\\small
           University of Copenhagen,             Laboratory of Actuarial Mathematics\\\small
            Universitetsparken 5,
            DK-2100 Copenhagen, Denmark\\\small
           \texttt{mikosch@math.ku.dk}
\and Olivier Wintenberger 
  \\\small       Centre De Recherche en Math\'ematiques de la D\'ecision
   UMR CNRS 7534
\\\small Universit\'e de Paris-Dauphine, Place du Mar\'echal De Lattre De Tassigny,
          \\\small  75775 Paris Cedex 16, France
         \\\small   \texttt{owintenb@ceremade.dauphine.fr}
            }

\date{}



\maketitle

\begin{abstract}
The aim of this paper is to provide conditions which ensure that the
affinely transformed partial sums of a strictly stationary process
converge in \ds\ to an infinite variance stable \ds . Conditions for
this \con\ to hold are known in the literature. However, most of
these results are qualitative in the sense that the parameters of
the limit \ds\ are expressed  in terms of some limiting \pp . In
this paper we will be able to determine the parameters of the
limiting stable \ds\ in terms of some tail characteristics of the
underlying stationary \seq . We will apply our results to some
standard \ts\ models, including the \garch\ process and its squares,
the \sv\ models and  solutions to \sre s.
\end{abstract}



\section{Introduction}\label{sec:1}
\setcounter{equation}{0} Whereas there exists a vast amount of
papers and books on the limit theory for sums $S_n=X_1+\cdots +X_n$
of finite variance strictly stationary \seq s $(X_t)$, less
attention has been given to the case of sums of infinite variance
stationary \seq  s. Following classical work (for example, Gnedenko
and Kolmogorov \cite{gnedenko:kolmogorov:1954}, Feller
\cite{feller:1971}, Petrov \cite{petrov:1995}), we know that an iid
\seq\ $(X_t)$ satisfies the limit relation \beqq\label{eq:2}
a_n^{-1} (S_n-b_n)\std Y_\alpha\,, \eeqq for suitable constants
$a_n>0$, $b_n\in \bbr$ and an infinite variance $\alpha$-stable \rv\
$Y_\alpha$ \fif\ the \rv\ $X=X_1$ has a \ds\ with \regvary\ tails
with index { $-\alpha\in (-2,0)$,} i.e., there exist constants $p,q\ge 0$
with $p+q=1$ and a \slvary\ \fct\ $L$  \st \beqq\label{eq:1}
\P(X>x){ \sim} p\,\dfrac{L(x)}{x^\alpha}\,\quad\mbox{and}\quad  \P(X\le
-x)\sim q\,\dfrac{L(x)}{x^\alpha}\,,\quad \xto \,. \eeqq This
relation is often referred to as {\em tail balance
  condition.} It will be convenient to refer to $X$ and its \ds\ as
{\em \regvary\ with index~$\alpha$}.

The limit relation \eqref{eq:2} is a benchmark result for weakly
dependent stationary \seq s with \regvary\ marginal \ds . However,
in the presence of dependence, conditions for the  \con\ of the
partial sums towards a stable limit are in general difficult to
obtain, unless some special structure is assumed. Early on,
$\alpha$-stable limit theory has been established for the partial
sums of linear processes $(X_t)$ with iid \regvary\ noise with index
$\alpha\in (0,2)$. Then the  linear process $(X_t)$ has \regvary\
marginals, each partial sum $S_d$, $d\ge 1$, is \regvary\ with index
$\alpha$
 and $(S_n)$ satisfies \eqref{eq:2} for
suitable $(a_n)$ and $(b_n)$. These results, the corresponding limit
theory for the partial sums $S_n$ and the sample autocovariance
\fct\ of linear processes were proved in a series of papers by Davis
and Resnick
\cite{davis:resnick:1985,davis:resnick:1985a,davis:resnick:1986}.
They exploited the relations between \regvar\ and the weak \con\ of
the \pp es $N_n=\sum_{t=1}^n \vep_{a_n^{-1}X_t}$, { where
$\vep_x$ denotes Dirac \ms\ at $x$.} Starting from the
\con\ $N_n\std N$, they used  a continuous mapping argument acting
on the points of the processes $N_n$ and $N$ in conjunction with the
series \rep\ of infinite variance stable \rv s. Their proofs heavily
depend on the linear dependence structure. A different, not \pp\
oriented,  approach was chosen by Phillips and Solo
\cite{philips:solo:1992} who decomposed the partial sums of the
linear process into an iid sum part and a negligible remainder term.
Then the limit theory for the partial sums follows from the one for
iid \seq s with \regvary\ marginal \ds .  The first result on stable
limits for stationary processes more general than linear models,
assuming suitable conditions for non-Gaussian limits, was proved by
Davis \cite{davis:1983}. Davis's ideas were further developed for
mixing sequences by Denker and Jakubowski
\cite{denker:jakubowski:1989} and Jakubowski and Kobus
\cite{jakubowski:kobus:1989}. The latter paper provides a formula
for the stable limit for sums of stationary sequences which are
$m$-dependent and admit local clusters of big values. A paper by
Dabrowski and Jakubowski \cite{dabrowski:jakubowski:1994} opened yet
another direction of studies: stable limits for associated
sequences.

Results for special non-linear \ts\ models, exploiting the structure
of the model,  were proved later on.  Davis and Resnick
\cite{davis:resnick:1996} and Basrak et al.
\cite{basrak:davis:mikosch:1999} studied the sample autocovariances
of bilinear processes with heavy-tailed and light-tailed noise,
respectively. Mikosch and Straumann \cite{mikosch:straumann:2006}
proved limit results for sums of stationary martingale differences
of the form $X_t=G_t\,Z_t$, where $(Z_t)$ is an iid \seq\ with
\regvary\ $Z_t$'s with index $\alpha\in (0,2)$, $(G_t)$ is adapted
to the filtration generated by $(Z_s)_{s\le t}$ and
$\E|G_t|^{\alpha+\delta}<\infty$ for some $\delta>0$. Stable limit
theory for the sample autocovariances of solutions to stochastic
recurrence equations, GARCH processes and \sv\ models was considered
in Davis and Mikosch \cite{davis:mikosch:1998,davis:mikosch:2001},
Mikosch and \sta\ \cite{mikosch:starica:2000}, Basrak et al.
\cite{basrak:davis:mikosch:2002}; see the survey papers Davis and
Mikosch
\cite{davis:mikosch:2009a,davis:mikosch:2009b,davis:mikosch:2009c}.

The last mentioned results are again based on the weak \con\ of the
\pp es $N_n=\sum_{t=1}^n \vep_{a_n^{-1}X_t}$ in combination with
continuous mapping arguments. The results make heavy use of the fact
that any $\alpha$-stable \rv , $\alpha\in (0,2)$, has a series \rep
, involving the points of a Poisson process. A general \asy\ theory
for partial sums of strictly stationary processes, exploiting the
ideas of \pp\ \con\ mentioned above, was given in Davis and Hsing
\cite{davis:hsing:1995}. The conditions in Davis and Hsing
\cite{davis:hsing:1995} are relatively straightforward to verify for
various concrete models. However, the $\alpha$-stable limits are
expressed as infinite series of the points of a Poisson process.
This fact makes it difficult to identify the parameters of the
$\alpha$-stable \ds s: these parameters are \fct s of the \ds\ of
the limiting \pp .

Jakubowski \cite{jakubowski:1993,jakubowski:1997} followed an
alternative approach based on classical blocking and mixing
techniques for partial sums of weakly dependent \rv s. A basic idea
of these papers consists of approximating the \ds\ of the sum
$a_n^{-1}S_n$ by the sum of the iid block sums
$(a_n^{-1}S_{mi})_{i=1,\ldots,k_n}$ \st\ $k_n=[n/m]\to \infty$ and
$S_{mi}\eqd S_m$. Then one can use the full power of classical
summation theory for row sums of iid  triangular arrays. It is also
possible to keep under control clustering of big values and
calculate the parameters of the $\alpha$-stable limit in terms of
quantities depending on the finite-dimensional distributions of the
underlying stationary process. Thus the direct method is in some
respects advantageous over the point process approach.

At a first glance, the conditions and results in Jakubowski
\cite{jakubowski:1993,jakubowski:1997} and Davis and Hsing
\cite{davis:hsing:1995} look rather different. Therefore we shortly
discuss these conditions in Section~\ref{sec:2} and argue that they
are actually rather close.
Our main result (Theorem~\ref{thm:clt}) is given in
Section~\ref{sec:main}. Using an argument going back to Jakubowski
\cite{jakubowski:1993,jakubowski:1997}, we provide an
$\alpha$-stable limit theorem for the partial sums of weakly
dependent infinite variance stationary \seq s. The proof only depends on the \chf s of
the converging partial sums. The result and its proof are new and
give insight into the dependence structure of a heavy-tailed
stationary \seq .  In Section~\ref{subsec:disc} we discuss the
conditions of Theorem~\ref{thm:clt} in detail. In particular, we
show that our result is easily applicable for strongly mixing \seq
s. In Section~\ref{sec:3} we explicitly calculate the parameters of
the $\alpha$-stable limits of the partial sums of the \garch\
process and its squares, solutions to \sre s, the \sv\ model and
symmetric $\alpha$-stable processes.

\section{A discussion of the conditions in $\alpha$-stable limit
  theorems}\label{sec:2}
\setcounter{equation}{0}
\subsection{Regular variation conditions}\label{subsec:rv}
We explained in Section~\ref{sec:1} that \regvar\ of $X$ with index
$\alpha\in (0,2)$ in the sense of \eqref{eq:1} is necessary and
sufficient for the limit relation \eqref{eq:2} with an
$\alpha$-stable limit $Y_\alpha$ for an iid \seq\ $(X_t)$. {The
necessity of \regvar\ of $X$ with index $\alpha\in (0,2)$ in the
case of dependent $X_i$'s is difficult to establish and, { in general,
incorrect; see Remark~\ref{rem:tm}.}
It is, however, natural to assume such a condition { as long as one takes the
conditions for an iid \seq\ as a benchmark result.}

Davis and Hsing  \cite{davis:hsing:1995} assume the stronger
condition that the strictly stationary \seq\ $(X_t)$ is {\em
\regvary\ with index $\alpha\in (0,2)$.} This means that the \fidi s
of $(X_t)$ have a jointly \regvary\ \ds\ in the following sense. For
every $d\ge 1$, there exists a non-null Radon \ms\ $\mu_d$ on the
Borel $\sigma$-field of $\ov \bbr^d\backslash\{\bf0\}$ (this means
that $\mu_d$ is finite on sets bounded away from zero), $\ov
\bbr=\bbr\cup\{\pm \infty\}$, \st\ \beqq\label{eq:5}
n\,\P(a_n^{-1}(X_1,\ldots,X_d)\in \cdot )\stv \mu_d(\cdot)\,, \eeqq
where $\stv$ denotes vague \con\ (see Kallenberg
\cite{kallenberg:1983}, Resnick \cite{resnick:1987}) and $(a_n)$
satisfies \beqq\label{eq:4} n\,\P(|X|>a_n)\sim 1\,. \eeqq The
limiting \ms\ has the property $\mu_d(xA)=x^{-\alpha}\mu_d(A)$,
$t>0$, for Borel sets $A$. We refer to $\alpha$ as the {\em index of
\regvar } of $(X_t)$ and its \fidi s. Note that Theorem 3 in
\cite{jakubowski:1994} provides conditions under which regular
variation of the one-dimensional marginals implies {\em joint}
regular variation (\ref{eq:5}).

Jakubowski \cite{jakubowski:1993,jakubowski:1997} does not directly
assume \regvar\ of $X$. However, his  condition U1 requires that the
normalizing \seq\ $(a_n)$ in \eqref{eq:2} is \regvary\ with index
$1/\alpha$. In \cite{jakubowski:1997} he also requires the
conditions ${\bf T}_+(d)$ and  ${\bf
  T}_-(d)$, $d\ge 1$, i.e., the existence of the limits
\beqq\label{eq:6} \lim_{\nto} n\,\P(S_d>a_n)=
b_+(d)\quad\mbox{and}\quad \lim_{\nto} n\,\P(S_d\le -
a_n)=b_-(d)\,,\quad d\ge 1\,. \eeqq If  $b_+(d)+b_-(d)>0$, the
\regvar\ of $(a_n)$ with index $1/\alpha$ is equivalent to \regvar\
of $S_d$ with index $\alpha$; see Bingham et al.
\cite{bingham:goldie:teugels:1987}. Condition U2 in
\cite{jakubowski:1997} restricts the class of all \regvary\ \ds s to
a subclass. The proof of Theorem~\ref{thm:clt} below shows that this
condition can be avoided. \bre\label{rem:2} Condition \eqref{eq:6}
is automatically satisfied for \regvary\ $(X_t)$, where
\beqq\label{cond:b} b_\pm (d)= \mu_d(\{\bfx\in \ov \bbr^d:
\pm(x_1+\cdots +x_d)>1\})\,. \eeqq Since $\mu_d$ is non-null for
every $d\ge 1$ and $\mu_d(tA)=t^{-\alpha}\mu(A)$, $t>0$, we have
$b_+(d)+b_-(d)>0$, $d\ge 1$. Since $(a_n)$ is \regvary\ with index
$1/\alpha$ it then follows that $S_d$ is \regvary\ with index
$\alpha$ for every $d\ge 1$. Since $(a_n)$ satisfies relation
\eqref{eq:4} it then follows that $b_+(1)=p$ and $b_-(1)=q$ with $p$
and $q$ defined in equation \eqref{eq:1}. In particular $p+q=1$. The
coefficients $b_+(d)$ and $b_-(d)$ for $d> 1$ can be considered as a
\ms\ of extremal dependence in the \seq\ $(X_t)$. The two benchmarks
are the iid case, $b_+(d)=p\,d$ and $b_-(d)=q\,d$ and the case
$X_i=X$ for all $i$, $b_+(d)=p\,d^\alpha$ and $b_-(d)=q\,d^\alpha$.
\ere Regular variation of a stationary \seq\ $(X_t)$ is a well
accepted concept in applied \pro y theory. One of the reasons for
this fact is that some of the important \ts\ models (ARMA with
\regvary\ noise, GARCH, solutions to \sre s, \sv\ models with
\regvary\ noise) have this property. Basrak and Segers
\cite{basrak:segers:2009} give some enlightening results about the
structure of \regvary\ \seq s.
In what follows, we will always assume: \\[1mm]
{\em Condition} {\bf (RV)}: The strictly stationary \seq\ $(X_t)$ is
\regvary\ with index $\alpha\in (0,2)$ in the sense of condition
\eqref{eq:5} with non-null Radon \ms s $\mu_d$, $d\ge 1$, and
$(a_n)$ chosen in \eqref{eq:4}.

\subsection{Mixing conditions}\label{subsec:22}
Assuming condition {\bf (RV)}, Davis and Hsing
\cite{davis:hsing:1995} require the mixing condition ${\mathcal
A}(a_n)$ defined in the following way. Consider the \pp\
$N_n=\sum_{t=1}^n\vep_{X_t/a_n}$  and assume that there exists a
\seq\ $m=m_n\to\infty$ \st\ $k_n=[n/m_n]\to\infty$, where $[x]$
denotes the integer part of $x$. The condition  ${\mathcal A}(a_n)$
requires that \beqq\label{eq:32} \E\ex^{-\int f\,dN_n} -
\left(\E\ex^{-\int f d N_m}\right)^{k_n}\to 0\,, \eeqq where $f$
belongs to a sufficiently rich class of non-negative measurable \fct
s on $\bbr$ \st\ the \con\ of the Laplace \fct al $\E\ex^{-\int
f\,dN_n}$ for all $f$ from this class ensures weak \con\ of $(N_n)$.
Relation \eqref{eq:32} ensures that $N_n$ can be approximated in law
by a sum of $k_n$ iid copies of~$N_m$, hence the weak limits of
$(N_n)$ must be infinitely divisible \pp es.

The condition ${\mathcal A}(a_n)$ is difficult to be checked
directly, but it follows from standard mixing conditions such as
strong mixing with a suitable rate. For future use, recall that the
stationary \seq\ $(X_t)$ is {\em strongly mixing with rate \fct }
$(\alpha_h)$ if \beao \sup_{A\in \sigma(\ldots,X_{-1},X_0)\,,B\in
  \sigma(X_h,X_{h+1},\ldots)}
\left|\P(A\cap B)-\P(A)\,\P(B)\right|=\alpha_h\to 0\,,\quad
h\to\infty\,. \eeao

Jakubowski \cite{jakubowski:1993} showed that \eqref{eq:2} with
$b_n=0$ and \regvary\ $(a_n)$   implies the condition
\beqq\label{eq:7} \max_{1\le k,l\le n\,,k+l\le n}\left| \E \ex^{ix
a_n^{-1}S_{k+l}} -
  \E\ex^{ix a_n^{-1}S_k}\E\ex^{ixa_n^{-1}S_l}\right|\to 0\,,\quad \nto ,\quad x\in\R
\eeqq which is satisfied for strongly mixing $(X_t)$. We also refer
to the discussion in Sections 4--6 of  \cite{jakubowski:1997} for
alternative ways of verifying \eqref{eq:7}. Under assumptions on the
\ds\ of $X$ more restrictive than \regvar\  it is shown that
\eqref{eq:2} implies the existence of a \seq\ $l_n\to\infty$ \st\
for any $k_n=o(l_n)$ the following relation holds \beqq\label{eq:35}
\left(\E\ex^{i\,x\,k_n^{-1/\alpha}(a_n^{-1}S_n)}\right)^{k_n} -
\E\ex^{i\,x\,a_n^{-1}S_n} \to 0\,,\quad x\in\bbr\,. \eeqq It is
similar to condition \eqref{eq:32} at the level of partial sums.

We will assume a similar mixing condition in terms of the \chf s of
the partial sums of $(X_t)$. Write \beao
\varphi_{nj}(x)=\E\ex^{ixa_n^{-1}S_j}\,,\quad j=1,2,\ldots\,,\quad
\varphi_n=\varphi_{nn}\,,\quad x\in \bbr\,. \eeao {\em Condition
{\bf (MX)}}. Assume that there exist $m=m_n\to\infty$ \st\
$k_n=[n/m]\to 0$ and \beqq\label{eq:mix}
\left|\varphi_n(x)-(\varphi_{nm}(x))^{k_n}\right|\to 0\,,\quad
\nto\,,\quad  x\in\bbr. \eeqq This condition is satisfied for a
strongly mixing \seq\ provided the rate \fct\ $(\alpha_h)$ decays
sufficiently fast; see Section~\ref{subsub:mx}. But \eqref{eq:mix}
is satisfied for classes of stationary processes much wider than
strongly mixing ones. Condition {\bf (MX)} is analogous to
${\mathcal A}(a_n)$. The latter condition is formulated in terms of
the Laplace \fct als of the underlying \pp es. It is motivated by
applications in \evt, where the weak \con\ of the \pp es is crucial
for proving limit results of the maxima and order statistics of the
samples $X_1,\ldots,X_n$. Condition {\bf (MX)} implies that the
partial sum processes  $(a_n^{-1}S_n)$ and
$(a_n^{-1}\sum_{i=1}^{k_n} S_{mi})$ have the same weak limits, where
$S_{mi}$, $i=1,\ldots,k_n$, are iid copies of $S_m$. This
observation opens the door to classical limit theory for partial
sums based on triangular arrays of independent \rv s. Since we are
dealing with the limit theory for the partial sum process
$(a_n^{-1}S_n)$ condition {\bf (MX)} is more natural than ${\mathcal
A}(a_n)$ which is only indirectly (via a non-trivial continuous
mapping argument acting on converging \pp es) responsible for the
\con\ of the normalized partial sum process  $(a_n^{-1}S_n)$.

\subsection{Anti-clustering conditions}\label{subsec:ac}
Assuming condition {\bf (RV)}, Davis and Hsing
\cite{davis:hsing:1995} require the anti-clustering condition
\beqq\label{eq:8} \lim_{d\to\infty}\limsup_{\nto} \P\Big( \max_{d\le
|i|\le m_n} |X_i|> x\,a_n\mid |X_0|>x\,a_n\Big)=0\,,\quad x>0\,,
\eeqq where, as before, $m=m_n\to\infty$ is the block size used in
the definition of the mixing condition ${\mathcal A}(a_n)$. It
follows from recent work by Basrak and Segers
\cite{basrak:segers:2009} that the index set $\{i:d\le |i|\le m_n\}$
can be replaced by $\{i:d\le i\le m_n\}$, reducing the efforts for
verifying \eqref{eq:8}. With this modification, a sufficient
condition for \eqref{eq:8} is then given by \beqq\label{eq:33}
\lim_{d\to\infty}\limsup_{\nto}n\, \sum_{i=d}^{m_n}\P\left(  |X_i|>
x\,a_n\,,|X_0|>x\,a_n\right)=0\,,\quad x>0\,. \eeqq Relation
\eqref{eq:33} is close to the anti-clustering condition $D'(x\,a_n)$
used in \evt ; see Leadbetter et al.
\cite{leadbetter:lindgren:rootzen:1983}, Leadbetter and Rootz\' en
\cite{leadbetter:rootzen:1988} and Embrechts et al.
\cite{embrechts:kluppelberg:mikosch:1997}, Chapter 5.

An alternative anti-clustering condition is (38) in Jakubowski
\cite{jakubowski:1997}: \beqq\label{eq:10}
\lim_{d\to\infty}\limsup_{x\to\infty}\limsup_{\nto} x^\alpha
\sum_{h=d}^{n-1} (n-h)\,\P(|X_0|>x\,a_n\,,|X_h|>x\,a_n)=0\,. \eeqq
Assuming \regvar\ of $X$ and defining $(a_n)$ as in  \eqref{eq:4},
we see that \eqref{eq:10} is implied by the condition \beao
\lim_{d\to\infty}\limsup_{x\to\infty}\limsup_{\nto}
n\,\sum_{h=d}^{n-1}\P(|X_h|>x\,a_n\,,|X_0|>x\,a_n)=0\,, \eeao which
is close to condition \eqref{eq:33}.

For our results we will need an anti-clustering condition as well.
It is hidden in assumption  \eqref{AC} in Theorem~\ref{thm:clt}; see
the discussion in Section~\ref{subsub:cond}.

\subsection{Vanishing small values conditions}\label{subsec:van}
Davis, Hsing, and Ja\-ku\-bows\-ki prove \con\ of the normalized
partial sums by showing that the limiting \ds\ is infinitely
divisible with a \levy\ triplet corresponding to an $\alpha$-stable
\ds . In particular, they need conditions to ensure that the sum of
the small values (summands) in the sum $a_n^{-1}S_n$ does not
contribute to the limit. Such a condition for a dependent \seq\
$(X_t)$ is often easily established  for $\alpha\in (0,1)$, whereas
the case $\alpha\in [1,2)$ requires some extra work.

Davis and Hsing \cite{davis:hsing:1995} assume the condition (3.2):
\beqq\label{eq:11} \lim_{\epsilon\to 0}\limsup_{\nto}
\P\left(\left|\sum_{t=1}^n
    X_t I_{\{|X_t|\le \epsilon a_n\}}- n\, \E XI_{\{|X|\le \epsilon
      a_n\}}\right| >x\,a_n\right)\, = 0, \quad x>0\,,
\eeqq for $\alpha\in (0,2)$. For an iid \seq\ $(X_t)$ the relation
\beao {\rm median}(a_n^{-1} S_n)-a_n^{-1} n \E [XI_{\{|X|\le
\epsilon a_n\}}]\to 0\,,\quad \epsilon>0\,, \eeao holds. Therefore
$a_n^{-1} n \E XI_{\{|X|\le \epsilon a_n\}}$ are the natural
centering constants for $a_n^{-1}S_n$ in stable limit theory. In the
case of dependent $(X_t)$, the choice of the latter centering
constants is less straightforward; it is dictated by truncation of
the points in the underlying weakly converging \pp es.

The analogous condition (35) in Jakubowski \cite{jakubowski:1997}
reads as follows: for each $x > 0$ \beqq\label{eq:12}
\lim_{\epsilon\downarrow 0}\limsup_{l\to\infty}\limsup_{\nto}
l^\alpha\, \P\Big(\Big|\sum_{t=1}^n [X_t\,I_{\{|X_t|\le \epsilon\,
      l\,a_n\}}-\E XI_{\{|X|\le
      \epsilon\,l\,a_n\}}]\Big|>x\,l\,a_n\Big) = 0.
\eeqq As shown in \cite{jakubowski:1997}, this condition is
automatically satisfied for $\alpha\in (0,1)$.

In our approach the anti-clustering condition {\bf (AC)} (see below)
is imposed on the sum of ``small" and ``moderate" values and we do
not need to verify conditions such as \eqref{eq:11} and
\eqref{eq:12}.

\section{Main result}\label{sec:main}
\setcounter{equation}{0} In this section we formulate and prove our
main result. Recall the \regvar\ condition {\bf (RV)} and the mixing
condition {\bf
  (MX)}
from Sections~\ref{subsec:rv} and \ref{subsec:22}. We will use the
following notation for any \rv\ $Y$: \beao \ov Y=(Y\wedge
2)\vee(-2)\,. \eeao Notice that $|\ov Y| = |Y|\wedge 2$ is
subadditive. \bth\label{thm:clt} Assume that $(X_t)$ is a strictly
stationary process satisfying the following conditions.
\begin{enumerate}
\item
The \regvar\ condition {\bf (RV)} holds for some $\alpha\in (0,2)$.
\item
The mixing condition {\bf (MX)} holds.
\item
The anti-clustering condition
\begin{equation}\label{AC}\tag{\bf AC}
\lim_{d\to\infty}\limsup_{n\to\infty} \dfrac n
m\sum_{j=d+1}^{m}\E\left|\overline{
    x\,a_n^{-1}(S_{j}-S_d)}\;\overline{ x\,a_n^{-1}
    X_{1}}\right|=0\,,\quad x\in\bbr\,,
\end{equation}
holds, where $m=m_n$ is the same as in {\bf (MX)}.
\item
The limits
\begin{equation}\label{TB}
\tag{\bf TB} \lim_{d\to \infty}(b_+(d)-b_+(d-1))=c_+\mbox{ and
}\lim_{d\to \infty}(b_-(d)-b_-(d-1))=c_-\,,
\end{equation}
exist. Here $b_+(d),b_-(d)$ are the tail balance parameters given in
\eqref{eq:6}.
\item
For  $\alpha>1$ assume $\E(X_1)=0$ and for $\alpha=1$,
\begin{equation}\label{CT}\tag{\bf CT}
\lim_{d\to\infty}\limsup_{n\to\infty}n\,|\E(\sin(a_n^{-1}S_d))|=0.
\end{equation}
\end{enumerate}
Then $c_+$ and $c_-$ are non-negative and $(a_n^{-1}S_n)$ converges
in \ds\ to an $\alpha$-stable \rv\ (possibly zero) with \chf\
$\psi_{\alpha}(x) = \exp( -|x|^{\alpha} \chi_{\alpha}(x, c_+,
c_-))$, where for $\alpha \neq 1$  the function $\chi_\alpha(x,c_+,
c_-), x \in \bbr$, is given by the formula
\[\dfrac{\Gamma(2-\alpha)}{1-\alpha}\,\Big((c_++c_-)\,\cos(\pi
\alpha/2)-i\,\sign(x) (c_+-c_-)\, \sin(\pi\, \alpha/2)\Big)\,,\]
while for $\alpha = 1$ one has
\[\chi_1(x,c_+, c_-) =
0.5\,\pi (c_++c_-) +i\,\sign(x)\,(c_+-c_-) \log |x|,\ x\in \bbr.
\]
\ethe We discuss the conditions of Theorem~\ref{thm:clt} in
Section~\ref{subsec:disc}. In particular, we compare them with the
conditions in Jakubowski \cite{jakubowski:1993,jakubowski:1997} and
Davis and Hsing~\cite{davis:hsing:1995}.
{ If the \seq\ $(X_n)$ is $m_0$-dependent for some integer
$m_0\ge 1$, i.e., the
$\sigma$-fields
$\sigma(\ldots,X_{-1},X_0)$ and $\sigma(X_{m_0+1},X_{m_0+2},\ldots)$ are
independent,
the conditions {\bf (MX), (AC)} and {\bf (TB)} of
Theorem~\ref{thm:clt} are automatic; see Section~\ref{subsec:md}.
The surprising fact that $c_+$ and $c_-$ are non-negative is explained
at the end of Section~\ref{subsec:tb}.}
\par
\bre\label{rem:tm}{ Although Theorem~\ref{thm:clt} covers a wide range of
strictly stationary \seq s (see in particular Section~\ref{sec:3})
condition {\bf (RV)} limits the applications to infinite variance
(in particular unbounded) \rv s $X_n$. The referee of this paper
pointed out the surprising fact that there exist strictly stationary
Markov chains $(Y_n)$, suitable bounded \fct s $f$ and a \seq\ $(a_n)$
with
$a_n=n^{1/\alpha} \ell(n)$ for some \slvary\ \fct\ $\ell$ \st\ the \seq\
of the normalized partial sums $(a_n^{-1}S_n)$ of the \seq\
$(X_n)=(f(Y_n))$ converges in \ds\ to an infinite
variance stable \rv . Then  {\bf (RV)} is obviously violated.
Such an example is contained in Gou\"ezel \cite{gouezel:2004},
Theorem 1.3.

Other examples of stable limits for sums of bounded stationary random variables
(of different nature - non-Markov and involving long-range dependence)
are given in \cite{surgailis:2004}, Theorems 2.1 (ii) and 2.2 (ii).
}\ere
\par
\bre\label{rem:eps}
It might be instructive to realize that in limit theorems for weakly dependent sequences properties of finite
dimensional distributions can be as bad as possible. For example it is very easy to build a $1$-dependent 
sequence having no moment of any order and such that its (centered and normalized) partial sums still converge 
to a stable law of order $\alpha \in (0,2]$. Let
\[ X_n=Y_n + \varepsilon_n - \varepsilon_{n-1},\]
where $(Y_n)$ is an iid sequence of $\alpha$-stable random variables and $(\varepsilon_n)$ is an iid sequence 
without any  moments (that is $E(|\varepsilon_0|^a)=+\infty$ for any positive $a$) and 
the two sequences are independent. Then the (centered and normalized) partial sums of $(X_n)$) have 
the same limit behavior (in distribution) as that of $(Y_n)$.
\ere

\subsection{Proof of Theorem~\ref{thm:clt}}\label{subsec:proof}
For any strictly stationary \seq\  $(X_t)_{t\in\Z}$ it will be
convenient to write \beao S_0=0\,,\quad S_n=X_1+\cdots +X_n\,,\quad
S_{-n}=X_{-n}+\cdots + X_{-1}\,,\quad n\ge 1\,. \eeao Let $S_{mi}$,
$i=1,2,\ldots,$ be iid copies of $S_m$. In view of \eqref{eq:mix}
the theorem is proved if we can show that
$(a_n^{-1}\sum_{i=1}^{k_n}S_{mi})$ has an $\alpha$-stable limit with
\chf\ $\psi_\alpha$. For such a triangular array, it is implied by the
relation
\begin{equation}\label{ta}
k_n(\varphi_{nm}(x)-1)\to \log\psi_\alpha(x)\,,\quad x\in\bbr\,.
\end{equation}
Indeed, notice first that the triangular array
$(a_n^{-1}S_{mi})_{i=1,\ldots,k_n}$ of iid \rv s satisfies the
infinite smallness condition. Then apply Lemma 3.5 in Petrov
\cite{petrov:1995} saying that for all $x\in\R$ and sufficiently
large $n$,
$$
\log
\varphi_{nm}(x)=\varphi_{nm}(x)-1+\theta_{nm}(\varphi_{nm}(x)-1)^2\,,
$$
where $|\theta_{nm}|\le 1$. Thus
$$
\left|k_n(\varphi_{nm}(x)-1)-k_n
 \log\varphi_{nm}(x)\right|\le k_n\,|\varphi_{nm}(x)-1|^2\le c\,k_n^{-1}\to 0\,.
$$
{\em Here and in what follows, $c$ denotes any positive constants.}
Our next goal is to find a suitable approximation to the \lhs\ in
\eqref{ta}.
\ble\label{ballou} Under {\bf (RV)} and {\bf (AC)} the following
relation holds:
\beqq\label{eq:ll}
\lim_{d\to\infty}\limsup_{n\to\infty}
\Big|k_n(\varphi_{nm}(x)-1)-n\,(\varphi_{nd}(x)-\varphi_{n,d-1}(x))\Big|=0\,,\quad
x\in\bbr\,. \eeqq
Moreover, if $(X_n)$ is $m_0$-dependent for some integer $m_0\ge 1$,
then
\beqq\label{eq:lla}
\lim_{n\to\infty}
\Big|k_n(\varphi_{nm}(x)-1)-n\,(\varphi_{nd}(x)-\varphi_{n,d-1}(x))\Big|=0\,,\quad
x\in\bbr\,,\quad d>m_0\,. \eeqq
\ele\noindent
The proof is given at the end of
this  section. By virtue of {\bf (RV)} and \eqref{eq:6}, $S_d$ is
regularly varying with index $\alpha\in (0,2)$; see
Remark~\ref{rem:2}. Therefore it belongs to the domain of attraction
of an $\alpha$-stable  law. Theorem 3 in Section XVII.5 of Feller
\cite{feller:1971} yields that for every $d\ge 1$ there exists an
$\alpha$-stable \rv\ $Z_\alpha(d)$ \st \beqq\label{eq:r}
a_n^{-1}\sum_{i=1}^n(S_{di}-e_{nd})\std Z_\alpha(d)\,\mbox{ where }\,e_{nd}=\begin{cases}0&\quad\alpha\neq 1,\\
\E(\sin(S_d/a_n))&\quad\alpha=1.\end{cases} \eeqq The limiting
variable $Z_\alpha(d)$ has the characteristic function
\[\widetilde{\psi}_{\alpha,d}(x) = \exp(-|x|^{\alpha}\chi_{\alpha}(x,
b_+(d),b_-(d)))\,,\quad x\in\bbr\,.\]
Applying Theorem 1  in Section XVII.5 of Feller \cite{feller:1971},
we find the equivalent relation
$$
n\,\big(\varphi_{nd}(x)\ex^{-ie_{nd}x}-1\big)\to \log \wt
\psi_{\alpha,d}(x)\,,\quad x\in\bbr\,,
$$
and exploiting condition \eqref{TB}, for $x\in\bbr$,
\beam
\lefteqn{n\,(\varphi_{nd}(x)\ex^{-ie_{nd}x}-\varphi_{n,d-1}(x)\ex^{-ie_{n,d-1}x})}\nonumber\\[2mm]
&\to& \log \wt \psi_{\alpha,d}(x)-\log \wt \psi_{\alpha,d-1}(x)\quad\mbox{as $\nto$}\label{h1}\\[2mm]
&\to& \log \psi_\alpha(x)\quad \mbox{as $d\to\infty$}\,.\nonumber
\eeam For $\alpha\neq1$ we have $e_{nd}=0$. Therefore \eqref{eq:ll}
implies
$$
k_n\,(\varphi_{nm}(x)-1)\to \log\psi_\alpha(x)\,,\quad x\in\bbr\,.
$$
This finishes the proof in this case. For $\alpha=1$ we  use the
same arguments but we have to take into account that $e_{nd}$ does
not necessarily vanish.
However, we have \beao
|\varphi_{nd}(x)-\varphi_{nd}(x)\ex^{-ie_{nd}x}|\le|1-\ex^{-ie_{nd}x}|\le
t\,  |e_{nd}|\,, \eeao and using \eqref{CT}, we obtain {\small \beao
\lim_{d\to\infty}\limsup_{n\to\infty}n\,|(\varphi_{nd}(x)-\varphi_{n,d-1}(x))-(\varphi_{nd}(x)\ex^{-ie_{nd}x}-\varphi_{n,d-1}(x)\ex^{-ie_{n,d-1}x})|=0\,.
\eeao} This proves the theorem.~\hfill\halmos
\noindent {\em Proof of Lemma \ref{ballou}.}
Consider the following telescoping sum for any $n$, $m\le n$  and
$d<m$: \beao
\varphi_{nm}(x)-1=\varphi_{nd}(x)-1+\sum_{j=1}^{m-d}(\varphi_{n,d+j}(x)-\varphi_{n,d-1+j}(x))\,.
\eeao By stationarity of $(X_t)$ we also have{\small \beao
\lefteqn{m\,(\varphi_{nd}(x)-\varphi_{n,d-1}(x)) = }\\[2mm]&&d\,(\varphi_{nd}(x)-
\varphi_{n,d-1}(x))+\sum_{j=1}^{m-d}\Big[ \E\ex^{ixa_n^{-1}(
S_{-d-j}- S_{-j})}- \E\ex^{ixa_n^{-1}( S_{-d-j+1}- S_{-j})}\Big]\,.
\eeao} Taking the difference between the previous two identities, we
obtain {\small \beao \lefteqn{
(\varphi_{nm}(x)-1)-m\,(\varphi_{nd}(x)-
   \varphi_{n,d-1}(x)) = }\\[2mm]
&&-(d-1)\,(\varphi_{n,d}(x)-1)+
d\,( \varphi_{n,d-1}(x)-1)+\\[2mm]
&& \sum_{j=1}^{m-d}\Big[ \varphi_{n,d+j}(x)-\E\ex^{ixa_n^{-1}(
S_{-d-j}- S_{-j})}- \varphi_{n,d-1+j}(x)+ \E\ex^{ixa_n^{-1}(
S_{-d-j+1}-S_{-j})}\Big]\,. \eeao} By stationarity, for any $k\ge1$,
$\varphi_{nk}(x)=\E\ex^{ixa_n^{-1} S_{-k}}$. Therefore any summand
in the latter sum can be written in the following form \beao
\lefteqn{ \E\Big(\ex^{ixa_n^{-1} S_{-d-j}}-\ex^{ixa_n^{-1}(
S_{-d-j}- S_{-j})}
-\ex^{ixa_n^{-1} S_{-d-j+1}}+\ex^{ixa_n^{-1}( S_{-d-j+1}- S_{-j})}\Big)}\qquad\qquad\qquad\\
&=&\E\Big( \ex^{ixa_n^{-1}
S_{-d-j}}\big(1-\ex^{-ixa_n^{-1}S_{-j}}\big)
\,\big(1-\ex^{-ixa_n^{-1}X_{-d-j}}\big)\Big)\,. \eeao Using the fact
that $x\to\exp(ix)$ is a $1$-Lipschitz function bounded by 1, the
absolute value of the expression on the \rhs\ is bounded by
$$
\E\big((|xa_n^{-1}S_{-j}|\wedge 2)\,(|xa_n^{-1}X_{-d-j}|\wedge
2)\big)\,.
$$
Collecting the above identities and bounds, we finally arrive at the
inequality
\beao
\lefteqn{|k_n\,(\varphi_{nm}(x)-1)-n\,(\varphi_{nd}(x)-\varphi_{n,d-1}(x))|}\qquad\qquad\qquad\nonumber\\[2mm]
&\le& k_n\,\,(d-1)\,| \varphi_{nd}(x)-1| + k_n\,d\,|\varphi_{n,d-1}(x)-1|\\[2mm]
&& + k_n\,\sum_{j=d+1}^{m}
\E\left|\ov {xa_n^{-1}(S_{j}-S_d)}\;\ov{xa_n^{-1} X_{1}}\right|.\nonumber\\
\eeao The last term on the \rhs\ converges to zero in view of
assumption \eqref{AC} when first $\nto$ and then $d\to\infty$.
In the $m_0$-dependent case, the last term on the \rhs\ converges to
zero whenever $d>m_0$ and $\nto$.
To prove that the first two terms also
converge to zero, let us notice that, under {\bf (RV)}, $n (\varphi_{nd}(x) -1 ) \to
\chi_{\alpha}(x,b_+(d),b_-(d))$ and so
\beao \lim_{n\to\infty} k_n (d-1) (\varphi_{nd}(x) - 1) = \lim_{n\to\infty}
{m_n}^{-1}(d-1)n (\varphi_{nd}(x) -1 ) = 0.\eeao This proves the lemma.\hfill\halmos

\bre Balan and Louhichi \cite{balan:louhichi:2008} have taken a
similar approach to prove limit theorems for triangular arrays of
stationary \seq s with infinitely divisible limits. Their paper
combines ideas from Jakubowski \cite{jakubowski:1997}, in particular
condition \eqref{TB}, and the \pp\ approach in Davis and Hsing
\cite{davis:hsing:1995}. They work under a mixing condition close to
${\mathcal A}(a_n)$. One of their key results (Theorem 2.6) is the
analog of Lemma~\ref{ballou} above. It is formulated in terms of the
Laplace \fct als of \pp es instead of the \chf s of the partial
sums. Then they sum the points in the converging \pp es and in the
limiting \pp\ to get an infinitely divisible limit. The sum of the
points of the limiting process represent an infinitely divisible
\rv\ by virtue of the \levy -It\^o \rep . As in Davis and Hsing
\cite{davis:hsing:1995} the method of proof is indirect, i.e., one
does not directly deal with the partial sums, and therefore the
results are less explicit. \ere

\subsection{A discussion of the conditions of
  Theorem~\ref{thm:clt}}\label{subsec:disc}

\subsubsection{Condition $5$}
It is a natural centering condition for the normalized partial sums
in the cases $\alpha=1$ and $\alpha\in (1,2)$. In the latter case,
$\E|X|<\infty$, and therefore $\E X=0$ can be assumed without loss
of generality. As usual in stable limit theory, the case $\alpha=1$
is special and therefore we need condition \eqref{CT}. It is
satisfied if $S_d$ is symmetric for every $d$.

\subsubsection{Condition \eqref{TB}}\label{subsec:tb}
If $c_++c_-=0$ the limiting stable \rv\ is zero.
 For example, assume $X_n=Y_n-Y_{n-1}$ for an iid  \regvary\
\seq\ $(Y_n)$ with index $\alpha\in (0,2)$. Then $S_d=Y_d-Y_0$ is
symmetric and \regvary\ with index $\alpha$. By the definition of
$(a_n)$, $b_+(d)=b_-(d)=0.5$, hence $c_+=c_-=0$. Of course,
$a_n^{-1}S_n\stp 0$.

In the context of Theorem~3.1 in Jakubowski \cite{jakubowski:1997}
(although the conditions of that result are more restrictive as
regards the tail of $X$) it is shown that \eqref{TB} is {\em
necessary} for \con\ of $(a_n^{-1}S_n)$ towards a stable limit.
Condition \eqref{TB} can be verified for various standard \ts\
models; see Section~\ref{sec:3}. The meaning of this condition is
manifested in Lemma~\ref{ballou}. It provides the link between the
\regvar\ of the \rv s $S_d$ for every $d\ge 1$ (this is a property
of the \fidi s  of the partial sum process $(S_d)$) and the \levy\
\ms\ $\nu_\alpha$ of the $\alpha$-stable limit. Indeed, notice that
\eqref{TB} implies that, for every $x>0$, with $b_+(0)=0$, \beao
c_+\,x^{-\alpha}&=& \lim_{d\to\infty} \frac{b_{+}(d)}{d}\,x^{-\alpha}\\[2mm]
&=&\lim_{d\to\infty}
\frac 1 d \sum_{i=1}^{d} (b_+(i)- b_+(i-1))\,x^{-\alpha}\\[2mm]&=&
\lim_{d\to\infty}\frac 1d \sum_{i=1}^d
\lim_{\nto}n\,(\P(S_i>x\,a_n)-\P(S_{i-1}>x\,a_n))\,, \eeao and a
similar relation applies to $c_- x^{-\alpha}$. Then \beao
\nu_\alpha(x,\infty)=c_+\,x^{-\alpha}\quad \mbox{and}\quad
\nu_\alpha(-\infty,-x)=c_-\,x^{-\alpha}\,,\quad x>0\,, \eeao
determine the \levy\ \ms\ $\nu_\alpha$ of the $\alpha$-stable limit
\ds\ with the characteristic function $\psi_{\alpha}$ given in
Theorem \ref{thm:clt}. In particular, Lemma~\ref{ballou} implies
that as $n\to\infty$ \beqq\label{eq:ld}
\left.\begin{array}{rl} k_n\,\P(S_m>x\,a_n)\to &\nu_\alpha(x,\infty),\\
k_n\,\P(S_m\le -x\,a_n)\to &\nu_\alpha(-\infty,-x],
\end{array}\right.\,\quad x>0\,.
\eeqq The latter relation opens the door to the limit theory for
partial sums of triangular arrays of iid copies
$(S_{mi})_{i=1,\ldots,k_n}$, of $S_m$. Notice that the relations
\eqref{eq:ld} are of \ld s type in the sense of
\cite{jakubowski:1997}. We refer to
\cite{jakubowski:nagaev:zaigraiev:1997} for their multi-dimensional
counterparts.

Let us notice that although one cannot ensure that $b_\pm(d)\ge
b_\pm(d-1)\ge 0$ for sufficiently large $d$, the constants $c_+,c_-$
are non-negative. It is immediate from the observation that
\begin{eqnarray*}
c_{\pm} &=& \lim_{d\to\infty} (b_{\pm}(d)- b_{\pm}(d-1)) \\
&=& \lim_{d\to\infty} \frac 1 d \sum_{i=1}^{d} (b_{\pm}(i)-
b_{\pm}(i-1)) = \lim_{d\to\infty} \frac{b_{\pm}(d)}{d} \geq 0.
\end{eqnarray*}

 \bre Recall the two benchmark examples of Remark \ref{rem:2}. If
$(X_t)$ is an iid \seq\ \regvary\ with index $\alpha>0$ the limits
$c_+=p$ and $c_-=q$ always exist and conditions {\bf (MX)}, {\bf
(AC)} are automatically satisfied. Then,  under {\bf (CT)}, we
recover the classical limit results for partial sums with
$\alpha$-stable limit. On the other hand, if $X_i=X$ for all $i$,
then $c_+=c_-=0$ if $0<\alpha<1$, $c_+=p$ and $c_-=q$ if $\alpha=1$
and $c_+$ and $c_-$ are not defined otherwise. This observation is
in agreement with the fact that $a_n^{-1}S_n= n^{1-1/\alpha} \ell(n)
X$ for some \slvary\ \fct\ $\ell$. \ere
\subsubsection{Sufficient conditions for \eqref{AC}}\label{subsub:cond}
Condition \eqref{AC} is close to the anti-clustering conditions in
\cite{davis:hsing:1995,jakubowski:1997} discussed in
Section~\ref{subsec:ac}. In what follows, we give some sufficient
conditions for \eqref{AC}. These conditions are often simple to
verify. \ble\label{lem:1} Assume the conditions of
Theorem~\ref{thm:clt} and that $(X_t)$ is strongly mixing with rate
\fct\ $(\alpha_h)$. Moreover, assume that there exists a \seq\
$r_n\to \infty$ \st\ $r_n/m_n\to 0$, $n\alpha_{r_n}\to 0$  and one
of the following three conditions is satisfied.
\begin{equation}\label{eq:a}
\begin{array}{rl}
\lim_{d\to\infty}\limsup_{\nto} n\,
\left[\sum_{i=d+1}^{r_n}\P(|X_i|> a_n\,,|X_1|>\,a_n)
\right.\qquad & \,\\
+\left.\P\left(\left|\sum_{i=d+1}^{r_n}
  X_iI_{\{|X_i|\le
    a_n\}}\right|>a_n\,,|X_1|>a_n\right)\right]=&0\,.
\end{array}
\end{equation}
or \beqq\label{eq:b} \lim_{d\to\infty}\limsup_{\nto}
n\,\P(\max_{i=d+1,\ldots,r_n} |X_i|>a_n/r_n\,,|X_1|>a_n)=0\,. \eeqq
or \beqq\label{eq:c} \lim_{d\to\infty}\limsup_{\nto}
n\,\P(|S_{r_n}-S_d|>a_n\,,|X_1|>a_n)=0\,. \eeqq Then \eqref{AC}
holds. \ele

\begin{proof} Let us recall that the function $y\mapsto |\overline{y}|$ is subadditive.
We decompose the sum in  \eqref{AC} as follows. \beao
k_n\,\left(\sum_{j=d+1}^{r_n}+\sum_{j=r_n+1}^m\right)\E\left|\overline{
    x\,a_n^{-1}(S_{j}-S_d)}\;\overline{ x\,a_n^{-1}
    X_{1}}\right|=J_1(n)+J_2(n)\,.
\eeao We will deal with the two terms $J_1(n)$ and $J_2(n)$ in
different ways. For the sake of simplicity we assume $x=1$.

We start by bounding $J_2(n)$. \beao J_2(n)&\le &
k_n\,\sum_{j=r_n+1}^m\left[\ \E\left|\overline{
a_n^{-1}(S_{j}-S_{r_n})}\;\overline{ a_n^{-1}
    X_{1}}\right|+ \left|\overline{a_n^{-1}(S_{r_n}-S_d)}\right|\ \right]\\[2mm]
&=& k_n\,\sum_{j=r_n+1}^m\E\left|\overline{
a_n^{-1}(S_{j}-S_{r_n})}\;\overline{ a_n^{-1}
    X_{1}}\right|+ \,n\,
\E\left|\overline{ a_n^{-1}(S_{r_n}-S_d)}\;\overline{
    a_n^{-1} X_{1}}\right|\\[2mm]
&=&J_{21}(n)+J_{22}(n)\,. \eeao We bound a typical summand in
$J_{21}(n)$, using the strong mixing property
\begin{eqnarray*}
\lefteqn{\E\left|\overline{ a_n^{-1}(S_j-S_{r_n})}\;\overline{ a_n^{-1} X_{1}}\right|}\qquad\\
&=&\cov\left(|\overline{ a_n^{-1}(S_j-S_{r_n})}|,|\overline{
a_n^{-1} X_{1}}|\right)
+\E\left|\overline{ a_n^{-1}S_{j-r_n+1}}\right|\E\left|\overline{ a_n^{-1} X_{1}}\right|\\
&\le&c\,\alpha_{r_n}+\E\left|\overline{
a_n^{-1}S_{j-r_n+1}}\right|\E\left|\overline{ a_n^{-1}
X_{1}}\right|.
\end{eqnarray*}
Moreover, we have for $j>r_n$
\[
\left|\overline{a_n^{-1}(S_{j-r_n+1})}\right|\leq \sum_{i=1}^{j-r_n
+1} \left|\overline{a_n^{-1}X_i}\right|,\] hence
$$
\E\left|\overline{ a_n^{-1}S_{j-r_n+1}}\right|\E\left|\overline{
a_n^{-1} X_{1}}\right|\le j\left(\E\left|\overline{ a_n^{-1}
X_{1}}\right|\right)^2.
$$
Thus we arrive at the bound \beao J_{21}(n)&\le&
c\,n\alpha_{r_n}+c\,n\,m\,\left(\E\left|\overline{ a_n^{-1}
X_{1}}\right|\right)^2\,. \eeao Observe that \beao \E|\overline{
a_n^{-1} X_{1}}|&=& \E \big(a_n^{-1} |X_{1}|I_{\{a_n^{-1} |X_1|\le
2\}}\wedge
2 I_{\{a_n^{-1} |X_1|>2\}}\big)\\[2mm]
&\le &2 \P(|X_1|>2\,a_n)\,. \eeao Therefore, by definition of
$(a_n)$ and since $n\alpha_{r_n}\to 0$ by assumption, \beao
J_{21}(n)=O(n\alpha_{r_n})+ O(m/n)=o(1)\,. \eeao We also have \beao
J_{22}(n)&\le  &c\,n\,\P(|S_{r_n}-S_d|>2\,a_n\,,|X_1|>2\,a_n)\\[2mm]
&\le &c\,n\,\sum_{i=d+1}^{r_n}\P(|X_i|> a_n\,,|X_1|>2\,a_n)\\[2mm]
&&+c\,n\,\P\left(\left|\sum_{i=d+1}^{r_n}
  X_iI_{\{|X_i|\le
    a_n\}}\right|>2\,a_n\,,|X_1|>2\,a_n\right)\,.
\eeao and \beao J_{22}(n)&\le &c\, \P(\max_{i=d+1,\ldots,r_n}
|X_i|>2\,a_n/r_n\mid |X_1|>2\,a_n)\,. \eeao Thus, under any of the
assumptions \eqref{eq:a}--\eqref{eq:c},
$\lim_{d\to\infty}\limsup_{\nto}J_{22}(n)=0$. Finally, \beao
J_1(n)&\le & c\,k_n\,\sum_{j=d+1}^{r_n}
\P(|S_j-S_d|>2a_n\,,|X_1|>2a_n)\\[2mm]
&\le & \dfrac{r_n}{m} \,n\,\P(|X_1|>2a_n)=o(1)\,. \eeao Collecting
the bounds above we proved that \eqref{AC} holds.
\end{proof}
\subsubsection{Condition {\bf (MX)}}\label{subsub:mx}
As we have already discussed in Sections~\ref{subsec:rv} and
\ref{subsec:22}, the condition {\bf (MX)} is a  natural one in the
context of stable limit theory for dependent stationary \seq s.
Modifications of these conditions appear in Davis and Hsing
\cite{davis:hsing:1995}, Jakubowski
\cite{jakubowski:1993,jakubowski:1997}. We also discuss the
existence of \seq es $m=m_n\to\infty$ and $r=r_n\to\infty$ \st\
$r/m\to0$ and $m/n\to0$ to be used in Lemma~\ref{lem:1} in the context of strong mixing.
\ble\label{lem:conmx} Assume that $(X_t)$ is strongly mixing with
rate \fct\ $(\alpha_h)$. In addition, assume that there exists a
sequence $\epsilon_n\to 0$ satisfying
\begin{eqnarray}
\label{condalpha1} n\,\alpha_{[\epsilon_n(a_n^2/n\wedge n)]}\to 0\,.
\end{eqnarray}
Then {\bf (MX)} holds for some $m=m_n\to\infty$ with $k_n=[n/m]\to
\infty$. Moreover, writing $r_n=[\epsilon_n (a_n\wedge n)]$, then
$r_n/m_n\to 0$ for this choice of $(m_n)$ and
\beqq\label{condalpha2} n\,\alpha_{r_n}\to 0\,. \eeqq \ele\noindent
If the tail index $\alpha\le 1$,
then \eqref{condalpha2}
turns into $n\alpha_{[\epsilon_n n]}\to0$ which is  more
restrictive than \eqref{condalpha1}. On the other hand, if the tail
index $\alpha$ is close to $2$, \eqref{condalpha1} is not implied by
polynomial decay of the coefficients $\alpha_h$. Then a
subexponential decay condition of the type $\alpha_n\le
C\exp(-cn^b)$ for some $C,c,b>0$ implies \eqref{condalpha1}, and
then \eqref{condalpha2} follows.
\begin{proof}
We start by showing that \eqref{eq:mix} holds for a suitable \seq\
$(m_n)$. Let $\varphi_{nm\delta}$ be the \chf\ of
$a_n^{-1}\sum_{i=1}^{k_n}U_{m-\delta,i}$ for some $\delta=\delta_n$
and \beao U_{ji}=\sum_{k=(i-1)j+1}^{ij}X_k \eeao a block sum of size
$j$. Using that characteristic functions are Lipschitz functions
bounded by 1 and writing  $q=k_n m$, for $x\in\bbr$, \beao
|\varphi_q(x)-\varphi_{nm\delta}(x)|&\le& \E \Big(\Big|\frac
x{a_n}\sum_{j=1}^{k_n}U_{\delta j}\Big|\wedge
  2\Big)\\[2mm]
&\le& \E\left(\frac{x\delta n}{m a_n}|X_1|\wedge 2\right)\\[2mm]
&\le &\int_0^2\P(|X_1|>ma_n/(\delta n)s)ds\,. \eeao The \rhs\
approaches zero if $m_na_n/(\delta_n n)\to \infty$ as  $n\to
\infty$. Under this condition, the same arguments yield \beao
&&|(\varphi_{nm}(x))^{k_n}-(\varphi_{n,m-\delta}(x))^{k_n}|\to
0\quad\mbox{and}\quad |\varphi_q(x)-\varphi_{n}(x)|\to 0\,, \eeao as
soon as $a_n/m\to \infty$. Next we use a standard mixing argument to
bound \beao
\lefteqn{|\varphi_{nm\delta}(x)-(\varphi_{n,m-\delta}(x))^{[m/n]}|}\\[2mm]
&\le&|\varphi_{nm\delta}(x)-\varphi_{m-\delta}(x)\varphi_{n-m,m\delta}(x)|\\&&+
|\varphi_{m-\delta}(x)\varphi_{n-m,m\delta}(x)-(\varphi_{n,m-\delta}(x))^{[m/n]}|.
\eeao The first term on the \rhs\ is the  covariance of bounded
Lipschitz functions of $S_{m-\delta}$ and $S_n-S_{m}$. Hence it is
bounded by $\alpha_\delta$. Iterative use of this argument,
recursively on distinct blocks,  shows that the \rhs\ is of the
order $(n/m)\alpha_\delta$. Thus we proved that \eqref{eq:mix} is
satisfied if \beqq\label{eq:ra} n/m\alpha_\delta \to 0\,,\quad
ma_n/(\delta n)\to \infty\quad\mbox{and}\quad a_n/m\to \infty\,.
\eeqq Choose $m_n=[\sqrt{\epsilon_n}(a_n\wedge n)]$, $\delta=[m^2/n]$
and assume \eqref{condalpha1}. Then \eqref{eq:ra} holds,
\eqref{eq:mix} is satisfied and $m/n\sim \sqrt{\delta/n}\to0$.

Finally, if \eqref{condalpha1} is satisfied choose
$r_n=[\epsilon_n(a_n\wedge n)]$. Then $n\alpha_{r_n}\to 0$ and $r_n/m\to0$ are
automatic.
\end{proof}

\section{Examples}\label{sec:3}
  \setcounter{equation}{0}
\subsection{$m_0$-dependent \seq s}\label{subsec:md}
Consider a strictly stationary \seq\ $(X_n)$  satisfying
condition {\bf (RV)} and which is $m_0$-dependent for some integer
$m_0\ge 1$. In this case, $\alpha_h=0$ for $h>m_0$. Then, by virtue of
Lemma~\ref{lem:conmx} condition {\bf (MX)} is satisfied
for any choice of \seq s $(m_n)$  \st\ $m_n\to\infty$ and $m_n=o(n)$.
Moreover,
${\bf (AC)}$ follows from Lemma~\ref{lem:1} for any $(r_n)$
\st\ $r_n\to\infty$ and $r_n=o(m_n)$. We verify the
validity of condition \eqref{eq:b}.  Then for $(r_n)$ growing
sufficiently slowly,
\beao\lefteqn{
n\,\P(\max_{i=d+1,\ldots,r_n} |X_i|>a_n/r_n\,,|X_1|>a_n)}\\[2mm]&\le&
n\,r_n \,P(|X_1|>a_n)\,P(|X_1|>a_n/r_n)\\[2mm]
&=& O(r_n \,P(|X_1|>a_n/r_n))=o(1)\,.
\eeao
Thus, in the $m_0$-dependent case we have the following special case of
Theorem~\ref{thm:clt}.
\bpr\label{prop:a}
Assume that $(X_t)$ is a strictly stationary $m_0$-dependent \seq\ for
some $m_0\ge 1$ which also satisfies condition  {\bf (RV)} for some
$\alpha\in (0,2)$. Moreover, assume $\E(X_1)=0$ for $\alpha>1$ and
$X_1$ is symmetric for $\alpha=1$. Then
the conclusions of Theorem~\ref{thm:clt} hold with
$c_+=b_+(m_0+1)-b_+(m_0)$ and $c_-=b_-(m_0+1)-b_-(m_0)$.
\epr
\begin{proof}
We have already verified conditions {\bf (MX)} and {\bf (AC)} of
Theorem~\ref{thm:clt}.
Following the lines of the proof of Theorem~\ref{thm:clt} with
$e_{nd}=0$,
we arrive at
\eqref{h1}
for every $d\ge 1$.  In view of the second part of Lemma~\ref{ballou},
the \rhs\ of  \eqref{h1} is independent of $d$ for $d>m_0$ as the
limit of $k_n (\varphi_{nm}(x)-1)$ as $\nto$. This finishes the proof
by taking  $d=m_0+1$.
\end{proof}
We mention in passing that we may conclude from the  proof of
Proposition~\ref{prop:a} that condition {\bf (TB)} is satisfied
since $c_+=b_+(d+1)-b_+(d)$ and $c_-=b_-(d+1)-b_-(d)$ for $d>m_0$.
This is a fact which is not easily seen by direct calculation on the
tails of $S_d$, $d>m_0$.

\subsection{The \sv\ model}\label{subsec:31}
The \sv\ model is one of the  standard econometric models for
financial returns of the form \beao X_t=\sigma_t\,Z_t\,, \eeao where
the volatility \seq\ $(\sigma_t)$ is strictly stationary independent
of the iid noise \seq\ $(Z_t)$. See e.g. Andersen et al.
\cite{andersen:davis:kreiss:mikosch:2009} for a recent reference on
\sv\ models or the collection of papers \cite{shephard:2005}.
\subsubsection*{Conditions {\bf (RV)}, \eqref{TB} and \eqref{CT}}
We assume that $Z$ is \regvary\ with index $\alpha>0$, implying that
$(Z_t)$ is \regvary . We also assume that $\E\sigma^{p}<\infty$ for
some $p>\alpha$. Under these assumptions it is known (see Davis and
Mikosch \cite{davis:mikosch:2001}) that $(X_t)$ is \regvary\ with
index $\alpha$, and the limit \ms\ $\mu_d$ in \eqref{eq:5} is given
by \beqq\label{eq:16} \mu_d(dx_1,\ldots,dx_d)=\sum_{i=1}^d
\la_\alpha(dx_i) \prod_{i\ne j} \vep_0(dx_j)\,, \eeqq where $\vep_x$
is Dirac \ms\ at $x$, \beao \la_\alpha(x,\infty)=\wt
p\,x^{-\alpha}\quad\mbox{and}\quad \la_\alpha(-\infty,-x]=\wt q\,
x^{-\alpha}\,,\quad x>0\,, \eeao and \beqq\label{eq:tildep} \wt
p=\lim_{\xto} \dfrac{\P(Z>x)}{\P(|Z|>x)} \quad\mbox{and}\quad \wt
q=\lim_{\xto} \dfrac{\P(Z\le -x)}{\P(|Z|>x)}\,, \eeqq are the tail
balance parameters of $Z$. This means that the \ms s $\mu_d$ are
supported on the axes as if the \seq\ $(X_t)$ were iid \regvary\
with tail balance parameters $\wt p$ and $\wt q$. By virtue of
\eqref{eq:16} and \eqref{cond:b} we have $b_{+}(d)= \wt p\,d$ and
$b_-(d)=\wt q\,d$, hence $c_+= \wt p$ and $c_-=\wt q$. We also
assume $\E Z=0$ for $\alpha>1$. Then $\E X=0$. If $\alpha=1$ we
assume $Z$ symmetric. Then $S_d$ is symmetric for every $d\ge 1$ and
\eqref{CT} is satisfied.
\subsubsection*{Conditions {\bf (MX)} and \eqref{AC} }
In order to meet {\bf (MX)} we assume that $(\sigma_t)$ is strongly
mixing with rate \fct\ $(\alpha_h)$. It is well known (e.g. Doukhan
\cite{doukhan:1994}) that $(X_t)$ is then strongly mixing with rate
\fct\ $(4\alpha_h)$.

It is common use in financial econometrics to assume that $(\log
\sigma_t)$ is a Gaussian linear process. The mixing rates for
Gaussian linear processes are well studied. For example, if $(\log
\sigma_t)$ is a Gaussian ARMA process then $(\alpha_h)$ decays
exponentially fast. We will assume this condition in the sequel.
Then we may apply Lemma~\ref{lem:conmx} with $r_n=n^{\gamma_1}$,
$m_n=n^{\gamma_2}$, $0<\gamma_1<\gamma_2<1$ for sufficiently small
$\gamma_1$ and $\gamma_2$, to conclude that {\bf (MX)} holds and
$n\alpha_{r_n}\to 0$.

Next we verify \eqref{AC}. We have by Markov's inequality for small
$\epsilon>0$,
\begin{eqnarray*}
\lefteqn{n\,\P (\max_{i=d+1,\ldots, r_n} |X_i|> a_n/r_n,|X_1|>a_n)}\qquad \\[2mm]
&\le &
n\,\sum_{i=d+1}^{r_n}\P ( |X_i|>  a_n/r_n,|X_1|>a_n)\\[2mm]
&\le &n\,\sum_{i=d+1}^{r_n}\P ( \max(\sigma_i,\sigma_1) \min
(|Z_i|,|Z_1|)>a_n/r_n)\\[2mm]
&\le &n
\,(r_n/a_n)^{\alpha+\epsilon}\sum_{i=d+1}^{r_n}\E(\max(\sigma_i^{\alpha+\epsilon},\sigma_1^{
\alpha+\epsilon})) \\[2mm]
&\le & c\, n\,r_n^{1+\alpha +\epsilon}\, a_n^{-\alpha-\epsilon}\,.
\end{eqnarray*}
The \rhs\ converges to zero if we choose $\gamma_1$ and $\epsilon$
sufficiently small. This proves \eqref{eq:b} and by
Lemma~\ref{lem:1} also \eqref{AC}. \bpr\label{pr:1} Assume that
$(X_t)$ is a \sv\ model satisfying the following additional
conditions:
\begin{enumerate}
\item[\rm (a)]
$(Z_t)$ is iid \regvary\ with index $\alpha\in (0,2)$ and tail
balance parameters $\wt p$ and $\wt q$.
\item[\rm (b)] For $\alpha \in (1,2)$, $\E Z=0$, and for $\alpha=1$,
$Z$ is symmetric.
\item[\rm (c)]
$(\log \sigma_t)$ is a Gaussian ARMA process.
\end{enumerate}
Then the \sv\ process $(X_t)$ satisfies the conditions of
Theorem~\ref{thm:clt} with parameters $c_+=\wt p$ and $c_-=\wt q$
defined in \eqref{eq:tildep}. \epr Hence a \sv\ model with Gaussian
ARMA log-volatility \seq\ satisfies the same stable limit relation
as an iid \regvary\ \seq\ with index $\alpha\in (0,2)$ and  tail
balance parameters $\wt p$ and $\wt q$.

In applications it is common to study powers of the absolute values,
$(|X_t|^p)$, most often for $p=1,2$. We assume the conditions of
Proposition~\ref{pr:1}. Then the \seq\ $(X_t^2)$ is again a \sv\
process which is \regvary\ with index $\alpha/2\in (0,1)$. It is not
difficult to see that the conditions of Proposition~\ref{pr:1} are
satisfied for this \seq\ with $b_-(d)=0$ and $b_+(d)=d$, hence
$c_+=1$ and $c_-=0$.

A similar remark applies to $(|X_t|)$ with one exception: the
centering condition \eqref{CT} cannot be satisfied. This case
requires special treatment. However, the cases $\alpha\ne 1$ are
similar. For $\alpha<1$, $(|X_t|)$ is a \sv\ model satisfying all
conditions of Proposition~\ref{pr:1}. For $\alpha\in (1,2)$ we
observe that \beao a_n^{-1} \sum_{i=1}^n (|X_t|-\E
|X|)&=&a_n^{-1}\sum_{t=1}^n \sigma_t \,(|Z_t|-\E|Z|)
+a_n^{-1}\E|Z|\,\sum_{t=1}^n
(\sigma_t-\E\sigma)\\
&=&a_n^{-1}\sum_{t=1}^n \sigma_t \,(|Z_t|-\E|Z|) +o_\P(1)\,. \eeao
In the last step we applied the \clt\ to $(\sigma_t)$. Then the
process $(\sigma_t(|Z_t|-\E |Z|))$ is a \sv\ model satisfying the
conditions of  Proposition~\ref{pr:1} with $c_+=1$ and $c_-=0$.

\subsection{Solutions to \sre s}\label{subsec:32}
We consider the \sre\ \beqq\label{eq:20}
X_t=A_t\,X_{t-1}+B_t\,,\quad t\in\bbz\,, \eeqq where $((A_t,B_t))$
constitutes an iid \seq\ of non-negative \rv s $A_t$ and $B_t$.
Various econometric \ts\ models $(X_t)$ have this form, including
the squared ARCH(1) process and the volatility \seq\ of a \garch\
process; see Section~\ref{subsec:33}. The conditions $\E\log A<0$
and $\E|\log B|<\infty$ are sufficient for the existence of a
strictly stationary causal solution $(X_t)$ to \eqref{eq:20} { \st\
  $(X_n)_{n\le 0}$ and $((A_n,B_n))_{n\ge 1}$} are independent; see Kesten
\cite{kesten:1973}.
\subsubsection*{Condition {\bf (RV)}}
Kesten \cite{kesten:1973} and Goldie \cite{goldie:1991} showed under
general conditions that $X$ has almost precise power law tail in the
sense that \beqq\label{eq:22} \P(X>x)\sim c_0\,x^{-\alpha} \eeqq for
some constant $c_0>0$, where the value $\alpha$ is given by the
unique positive solution to the equation \beao \E A^\kappa=1\,,\quad
\kappa>0\,. \eeao We quote Theorem 4.1 in Goldie \cite{goldie:1991}
to get the exact conditions for \eqref{eq:22}.
\bth\label{thm:2}
Assume that $A$ is a non-negative \rv\ \st\ the conditional law of
$A$ given $A\ne 0$ is non-arithmetic and there exists $\alpha>0$
\st\
 $\E A^\alpha=1$, $\E (A^\alpha\log^+ A)<\infty$.
Then $-\infty \le \E\log A<0$ and $\E (A^\alpha \log A)\in
(0,\infty)$. Moreover, if $\E B^\alpha<\infty$, then a unique
strictly stationary causal solution $(X_t)$ to \eqref{eq:20} exists \st\
\eqref{eq:22} holds with constant
\beqq\label{eq:27} c_0=
\dfrac{\E[(B_1+A_1X_0)^\alpha-(A_1\,X_0)^\alpha] }{\alpha \, \E
(A^\alpha\log A)}\,. \eeqq \ethe\noindent
The condition of
non-arithmeticity of the \ds\ of $A$ is satisfied if $A$ has a
Lebesgue density. {\em In what follows, we assume that the
conditions of Theorem~\ref{thm:2} are satisfied.}

Iterating the defining equation \eqref{eq:20} and writing \beao
\Pi_t= A_1\cdots A_t\,,\quad t\ge 1\,, \eeao we see that
\beqq\label{eq:21} (X_1,\ldots,X_d)= X_0
\,(\Pi_1,\Pi_2,\ldots,\Pi_d)+R_d\,, \eeqq where $R_d$ is independent
of $X_0$. Under the assumptions of Theorem~\ref{thm:2}, the moments
 $\E A^\alpha$ and $\E B^\alpha$ are finite, hence $\E (R_d^\alpha)<\infty$ and
$\P(|R_d|>x)=o(\P(|X_0|>x))$. By a multivariate version of a result
of Breiman \cite{breiman:1965} (see Basrak et al.
\cite{basrak:davis:mikosch:2002}) it follows that the first term on
the \rhs\ of \eqref{eq:21} inherits the \regvar\ from $X_0$ with
index $\alpha$ and by a standard argument (see Jessen and Mikosch
\cite{jessen:mikosch:2006}, Lemma~3.12) it follows that
$(X_1,\ldots,X_d)$ and the first term on the \rhs\ of \eqref{eq:21}
have the same limit \ms\ $\mu_d$. Hence the \seq\ $(X_t)$ is
\regvary\ with index $\alpha$, i.e., condition {\bf (RV)} is
satisfied for $\alpha>0$ with $\E A^\alpha=1$.
\subsubsection*{Condition \eqref{TB}}
Next we want to determine the quantities $b_+(d)$. Choose $(a_n)$
\st\ $n\,\P(X>a_n)\sim 1$, i.e., $a_n=(c_0\,n)^{1/\alpha}$, and
write \beao T_d=\sum_{i=1}^d \Pi_i\,,\quad d\ge 1\,. \eeao We obtain
for every $d\ge 1$, by \eqref{eq:21}, \beao n\,\P(S_d> a_n)&\sim
&n\, \P( X_0\,T_d>a_n) \sim  n\,\P(X_0>a_n)\, \E (T_d^\alpha) \sim
\E (T_d^\alpha) =b_+(d)\,. \eeao Here we again used Breiman's result
\cite{breiman:1965} for $\P(X_0\,T_d>x)\sim \E (T_d^\alpha) \P
(X_0>x)$ in a modified form. In general, this result requires that
$\E (T_d^{\alpha+\delta})<\infty$ for some $\delta>0$. However, if
$\P(X>x)\sim c_0\,x^{-\alpha}$, Breiman's result is applicable under
the weaker condition $\E A^\alpha<\infty$; see Jessen and Mikosch
\cite{jessen:mikosch:2006}, Lemma 4.2(3). Of course, $b_-(d)=0$. We
mention that the values $b_\pm (d)$ do not change if $S_d$ is
centered by a constant.

Our next goal is to determine $c_+$. Since $\E A^\alpha=1$ we have
\beqq\label{eq:forget} b_+(d+1)-b_+(d)= \E[(1+T_d)^\alpha-
T_d^\alpha]\,. \eeqq The condition $\E A^\alpha=1$ and convexity of
the \fct\ $g(\kappa)=\E A^\kappa$, $\kappa>0$, imply that $\E\log
A<0$ and therefore \beao T_d\stas T_\infty = \sum_{i=1}^\infty
\Pi_i<\infty\,. \eeao Therefore the question arises as to whether
one may let $d\to\infty$ in \eqref{eq:forget} and replace $T_d$ in
the limit by $T_\infty$. This is indeed possible as the following
dominated \con\ argument shows.

If $\alpha\in (0,1]$ concavity of the \fct\ $f(x)=x^\alpha$ yields
that $(1+T_d)^\alpha- T_d^\alpha\le 1$ and then Lebesgue dominated
\con\ applies. If $\alpha\in (1,2)$, the mean value theorem yields
that \beao (1+T_d)^\alpha - T_d^\alpha=
\alpha\,(T_d+\xi)^{\alpha-1}\,, \eeao where $\xi \in (0,1)$. Hence
$(1+T_d)^\alpha - T_d^\alpha$ is dominated by the \fct\ $\alpha
[T_d^{\alpha-1}+1]$. By convexity of $g(\kappa)$, $\kappa>0$, we
have  $\E (A^{\alpha-1})<1$ and therefore \beao \E
(T_\infty^{\alpha-1})\le\sum _{i=1}^\infty \E(\Pi_i^{\alpha-1})=
\sum _{i=1}^\infty (\E( A^{\alpha-1}))^i= \E (A^{\alpha-1}) (1-\E
(A^{\alpha-1}))^{-1}<\infty\,. \eeao

An application of Lebesgue dominated \con\ yields for any $\alpha\in
(0,2)$ that \beqq\label{eq:25}
c_+=\lim_{d\to\infty}[b_+(d+1)-b_+(d)]= \E[(1+T_\infty)^\alpha-
T_\infty^\alpha]\in (0,\infty)\,. \eeqq \bre\label{rem:1} The
quantity $T_\infty$ has the stationary \ds\ of the solution to the
\sre \beao Y_t= A_t\,Y_{t-1}+1\,,\quad t\in \bbz\,. \eeao This
solution satisfies the conditions of Theorem~\ref{thm:2} and
therefore \beao \P(Y_0>x)=\P(T_\infty>x)\sim c_1\,x^{-\alpha}\,,
\eeao with constant \beao c_1=
\dfrac{\E[Y_1^\alpha-(A_1\,Y_0)^\alpha]}{\alpha\,\E (A^\alpha\log
A)} =\dfrac{\E[(1+A_1\,Y_0)^\alpha-(A_1\,Y_0)^\alpha]}{\alpha\,\E(
A^\alpha\log A)}\,. \eeao In particular, $\E(
T_\infty^\alpha)=\infty$. This is an interesting observation in view
of $c_+\in (0,\infty)$. It is also interesting to observe that the
limit relation \eqref{eq:25} implies that \beao \dfrac{b_+(d)}{d}=
\dfrac{\E (T_d^\alpha)}{d}= \E[d^{-1/\alpha} T_d]^\alpha \to
\E[(1+T_\infty)^\alpha-T_\infty^\alpha]\,, \eeao although
$d^{-1/\alpha}T_d\stas 0$. This relation yields some information
about the rate at which $T_d\stas T_\infty$. \ere
\subsubsection*{Condition {\bf (MX)}}
The stationary solution $(X_t)$ to the \sre\ \eqref{eq:20} is
strongly mixing with geometric rate provided that some additional
conditions are satisfied. For example, Basrak et al.
\cite{basrak:davis:mikosch:2002}, Theorem 2.8, assume that the \MC\
$(X_t)$ is  $\mu$-irreducible, allowing for the machinery for Feller
chains with drift conditions as for example explained in Feigin and
Tweedie \cite{feigin:tweedie:1985} or Meyn and Tweedie
\cite{meyn:tweedie:1993}. The drift condition can be verified if one
assumes that $A_t$ has polynomial structure; see Mokkadem
\cite{mokkadem:1990}. The latter conditions can be calculated for
GARCH and bilinear processes, assuming some positive Lebesgue
density for the noise in a neighborhood of the origin; see Basrak et
al. \cite{basrak:davis:mikosch:2002}, Straumann and Mikosch
\cite{straumann:mikosch:2006}.

{\em In what follows, we will assume that $(X_t)$ is strongly mixing
with geometric rate.} Then, by Lemma~\ref{lem:conmx}, we may assume
that we can choose $r_n=n^{\gamma_1}$, $m_n=n^{\gamma_2}$ for
sufficiently small values $0<\gamma_1<\gamma_2<1$. Then {\bf (MX)}
holds and $n\alpha_{r_n}\to 0$.
\subsubsection*{Condition \eqref{AC}}
We verify condition \eqref{eq:c} and apply Lemma~\ref{lem:1}. It
suffices to bound the quantities \beao
I_n(d)=\P(|S_{r_n}-S_d|>a_n\mid X_0>a_n)\,. \eeao Writing
$\Pi_{s,t}= \prod_{i=s}^t A_i$ for $s\le t$ and $\Pi_{st}=1$ for
$s>t$, we obtain \beao X_i= X_0\,\Pi_i + \sum_{l=1}^i
\Pi_{l+1,i}\,B_l= X_0\,\Pi_i +C_i\,,\quad i\ge 1\,. \eeao Then,
using the independence of $X_0$ and $C_i$, $i\ge 1$, applying
Markov's inequality for $\kappa<\alpha\wedge 1$ and Karamata's
theorem (see Bingham et al. \cite{bingham:goldie:teugels:1987}),
\beao I_n(d)&\le &\P\Big( X_0 \sum_{i=d+1}^{r_n}\Pi_i>a_n/2\mid
X_0>a_n\Big)
+\P\Big(\sum_{i=d+1}^{r_n}C_i>a_n/2\Big)\\[2mm]
&\le &c\,\dfrac{\E(X_0^\kappa
\,I_{\{X_0>a_n\}})}{a_n^\kappa\,\P(X_0>a_n)}
\sum_{i=d+1}^{r_n}E(\Pi_i^\kappa)
+c\,a_n^{-\kappa}\sum_{i=d+1}^{r_n}\sum_{l=1}^i (\E A^\kappa)^{i-l}\,\E B^\kappa\\[2mm]
&\le &c\,\sum_{i=d+1}^{\infty} (\E (A^\kappa))^i
+c\,r_n^{1+\kappa} a_n^{-\kappa}(1- \E A^\kappa)^{-1} \E B^\kappa\\[2mm]
&\le & c\,\big((\E (A^\kappa))^d +r_n^{1+\kappa} a_n^{-\kappa})\,.
\eeao Here we also used the fact that $\E(A^\kappa)<1$ by convexity
of the \fct\ $g(\kappa)=\E(A^\kappa)$, $\kappa>0$, and
$g(\alpha)=1$. Choosing $r_n=n^{\gamma_1}$ for $\gamma_1$
sufficiently small, we see that \beao
\lim_{d\to\infty}\limsup_{\nto}  I_n(d)=0\,. \eeao This proves
\eqref{eq:c}.
\subsubsection*{Condition $5$}
Since $X$ is non-negative, \eqref{CT} cannot be satisfied; the case
$\alpha=1$ needs special treatment. We focus on the case $\alpha\in
(1,2)$. It is not difficult to see that all calculations given above
remain valid if we replace $X_t$ by $X_t-\E X$, provided $\gamma_1$
in $r_n=n^{\gamma_1}$ is chosen sufficiently small.

We summarize our results. \bpr\label{pr:2} Under the conditions of
Theorem~\ref{thm:2} the \sre\ \eqref{eq:20} has a strictly
stationary solution $(X_t)$ which is \regvary\ with index $\alpha>0$
given by $\E[ A^\alpha]=1$. If $\alpha\in (0,1)\cup (1,2)$ and
$(X_t)$ is strongly mixing with geometric rate the conditions of
Theorem~\ref{thm:clt} are satisfied. In particular, \beao (c_0
n)^{-1/\alpha}\,(S_n-b_n)\std Z_\alpha\,,\quad \mbox{where}\quad
b_n=\left\{\barr{ll}
0& \mbox{for $\alpha\in (0,1)$}\,,\\[2mm]
\E S_n=n\E X &\mbox{for $\alpha\in (1,2)$}\,, \earr\right. \eeao
where the constant $c_0$ is given in \eqref{eq:27} and the
$\alpha$-stable \rv\  $Z_\alpha$ has \chf\ $\psi_{\alpha}(t) = \exp(
-|t|^{\alpha} \chi_{\alpha}(t, c_+, 0))$, where \beao c_+=
\E[(1+T_\infty)^\alpha-T_\infty^\alpha]\in (0,\infty)\quad
\mbox{and}\quad T_\infty=\sum_{i=1}^\infty A_1\cdots A_i\,. \eeao
\epr
\bre{ Analogs of of Proposition~\ref{pr:2} have recently been
  proved in Guivarc'h and Le Page \cite{guivarch:lepage:2008} in the one-dimensional case and in Buraczewski et
  al. \cite{buraczewski:damek:guivarch:2009}, Theorem 1.6, also in the
  multivariate case. The results are
  formulated for a non-stationary version of the process $(X_n)$
  starting at some fixed value $X_0=x$. (This detail is not essential
  for the limit theorem.)
The proofs are tailored
for the situation of \sre s and therefore different from those in this
paper where the proofs do not depend on some particular structure of
the underlying stationary \seq .
}
\ere
\subsection{ARCH(1) and GARCH(1,1) processes}\label{subsec:33}
In this section we consider the model \beao X_t=\sigma_t\,Z_t\,,
\eeao where $(Z_t)$ is an iid \seq\ with $\E Z=0$ and $\var(Z)=1$
and \beqq\label{eq:23} \sigma_t^2=\alpha_0+ (\alpha_1
Z_{t-1}^2+\beta_1)\,\sigma_{t-1}^2\,. \eeqq We assume that
$\alpha_0>0$ and the non-negative parameters $\alpha_1,\beta_1$ are
chosen \st\ a strictly stationary solution to the \sre\
\eqref{eq:23} exists, namely, \beqq\label{eq:00} -\infty \le \E\log
(\alpha_1\,Z^2+\beta_1)<0\,, \eeqq see Goldie \cite{goldie:1991},
cf. Mikosch and  \sta\ \cite{mikosch:starica:2000}. Then the process
$(X_t)$ is strictly stationary as well. It is called a GARCH(1,1)
process if $\alpha_1\beta_1>0$ and an  ARCH(1) process if
$\beta_1=0$ and $\alpha_1>0$. Notice that condition \eqref{eq:00}
implies that $\beta_1 \in [0,1)$.

As a matter of fact, these classes of processes fit nicely into the
class of \sre s considered in Section~\ref{subsec:32}. Indeed, the
squared volatility process $(\sigma_t^2)$ satisfies the \sre\
\eqref{eq:20} with $X_t=\sigma_t^2$, $A_t=\alpha_1
Z_{t-1}^2+\beta_1$ and $B_t=\alpha_0$.  Moreover, the squared
ARCH(1) process $(X_t^2)$ satisfies \eqref{eq:20} with $Y_t=X_t^2$,
$A_t=\alpha_1 Z_t^2$ and $B_t=\alpha_0\,Z_t^2$.

A combination of the  results in Davis and Mikosch
\cite{davis:mikosch:1998} for ARCH(1) and in Mikosch and \sta\
\cite{mikosch:starica:2000} for \garch\ with Proposition~\ref{pr:2}
above yields the following. \bpr\label{pr:3} Let $(X_t)$ be a
strictly stationary \garch\ process. Assume that $Z$ has a positive
density on $\bbr$  and that there exists $\alpha>0$ \st\
\beqq\label{eq:24} \E[(\alpha_1Z^2+\beta_1)^\alpha]=1\quad
\mbox{and}\quad \E[(\alpha_1 Z^2+\beta_1)^\alpha \log (\alpha_1
Z^2+\beta_1)]<\infty\,. \eeqq Then the  following statements hold.
\begin{enumerate}
\item[\rm (1)]
The stationary solution $(\sigma_t^2)$ to \eqref {eq:23} is
\regvary\ with index $\alpha$ and strongly mixing with geometric
rate. In particular, there exists a constant $c_1>0$, given in
\eqref{eq:27} with $A_1=\alpha_1 Z_0^2+\beta_1$ and $B_1=\alpha_0$
\st \beqq\label{eq:34} \P(\sigma^2>x)\sim c_1\,x^{-\alpha}\,. \eeqq
For $\beta_1=0$, the squared {\rm ARCH(1)} process $(X_t^2)$ is
\regvary\ with index $\alpha$ and strongly mixing with geometric
rate. In particular, there exists a constant $c_0>0$ given in
\eqref{eq:27} with $A_1=\alpha_1 Z_1^2$ and $B_1=\alpha_0Z_1^2$ \st
\beao \P(X^2>x)\sim c_0\,x^{-\alpha}\,. \eeao
\item[\rm (2)]
Assume $\alpha\in (0,1)\cup(1,2)$ in the \garch\ case. Then \beao
(c_1 n)^{-1/\alpha} \Big(\sum_{t=1}^n \sigma_t^2- b_n\Big)\std
Z_\alpha\,,
\quad \mbox{where}\quad b_n=\left\{\barr{ll}0& \alpha \in (0,1)\\[2mm]
 n\,\E(\sigma^2) & \alpha\in (1,2)\,,\earr\right.
\eeao and the $\alpha$-stable \rv\ $Z_\alpha$ has \chf\
$\psi_{\alpha}(t) = \exp( -|t|^{\alpha} \chi_{\alpha}(t, c_+, 0))$,
where \beao c_+=\E[(1+T_\infty)^\alpha-T_\infty^\alpha]\,\quad
\mbox{and}\quad T_\infty=\sum_{t=1}^\infty
\prod_{i=1}^t(\alpha_1\,Z_i^2+\beta_1)\,. \eeao
\item[\rm (3)]
Assume $\alpha\in (0,1)\cup(1,2)$ in the {\rm ARCH(1)} case. Then
\beao (c_0 n)^{-1/\alpha} \Big(\sum_{t=1}^n X_t^2- b_n\Big)\std \wt
Z_\alpha\,,
\quad \mbox{where}\quad b_n=\left\{\barr{ll}0& \alpha \in (0,1)\\[2mm]
 n\,\E(\sigma^2) & \alpha\in (1,2)\,,\earr\right.
\eeao and the $\alpha$-stable \rv\ $\wt Z_\alpha$ has \chf\
$\psi_{\alpha}(t) = \exp( -|t|^{\alpha} \chi_{\alpha}(t, c_+, 0))$,
where \beao c_+=\E[(1+\wt T_\infty)^\alpha-\wt
T_\infty^\alpha]\,\quad \mbox{and}\quad \wt
T_\infty=\sum_{t=1}^\infty \alpha_1^t\,\prod_{i=1}^t Z_i^2\,. \eeao
\end{enumerate}
\epr The  limit results above require that we know the constants
$c_1$ and $c_0$ appearing in the tails of $\sigma^2$ and $X^2$. For
example, in the ARCH(1) case, \beao c_0=
\dfrac{\E[(\alpha_0+\alpha_1\, X^2)^\alpha- (\alpha_1\,X^2)^\alpha]
\,\E|Z|^{2\alpha}}{\alpha\ \E[(\alpha_1 \,Z^2)^\alpha \log (\alpha_1
Z^2)]}\,. \eeao Moreover, (8.66) in  Embrechts et al.
\cite{embrechts:kluppelberg:mikosch:1997} yields that $ X^2\eqd
(\alpha_0/\alpha_1)\wt T_\infty\,. $ Hence the constant $c_0$ can be
written in the form \beao c_0= \dfrac{\alpha_0^\alpha\,\E[(1+\wt
T_\infty )^\alpha- \wt T_\infty^\alpha] \;\E|Z|^{2\alpha}}{\alpha\
\E[(\alpha_1 \,Z^2)^\alpha \log (\alpha_1 Z^2)]} =
\dfrac{\alpha_0^\alpha\,c_+\,\E|Z|^{2\alpha}}{\alpha\ \E[(\alpha_1
\,Z^2)^\alpha \log (\alpha_1 Z^2)]}\,. \eeao The moments of $|Z|$
can be evaluated by numerical methods given that one assumes that
$Z$ has a tractable Lebesgue density, such as the standard normal or
student densities. Using similar numerical techniques, the value
$\alpha$ can be derived from \eqref{eq:24}. The evaluation of the
quantity $\E[(\alpha_0+\alpha_1\, X^2)^\alpha-
(\alpha_1\,X^2)^\alpha]$ is a hard problem; Monte-Carlo simulation
of the ARCH(1) process is an option. In the general \garch\ case,
similar remarks apply to the constants $c_+$ and $c_1$ appearing in
the stable limits of the partial sum processes of $(\sigma_t^2)$.
Various other financial \ts\ models fit into the framework of \sre
s, such as the AGARCH and EGARCH models; see e.g. the treatment in
Straumann and Mikosch \cite{straumann:mikosch:2006} and the lecture
notes by Straumann \cite{straumann:2005}.

In what follows, we consider the \garch\ case and  prove stable
limits for the partial sums of $(X_t)$ and $(X_t^2)$.
\subsubsection*{Conditions {\bf (RV)},  {\bf (MX)} and \eqref{AC}}
{\em In what follows, we assume that $Z$ is symmetric, has a
positive Lebesgue density on $\bbr$ and there exists $\alpha>0$ \st\
\eqref{eq:24} holds.} Under these assumptions, it follows from
Mikosch and \sta\ \cite{mikosch:starica:2000} that $(X_t^2)$ is
\regvary\ with index $\alpha$  and strongly mixing with geometric
rate. By Breiman's \cite{breiman:1965} result we have in particular,
\beao \P(X^2>x)\sim \E|Z|^{2\alpha}\,\P(\sigma^2>x)\sim
\E|Z|^{2\alpha}\, c_1\,x^{-\alpha}\,. \eeao By definition of
multivariate \regvar , the \seq\ $(|X_t|)$ inherits \regvar\ with
index $2\alpha$ from $(X_t^2)$. By symmetry of $Z$  the \seq s
$(\sign (Z_t))$ and $(|Z_t|)$, hence $(\sign (X_t))$ and $(|X_t|)$,
are independent. Then an application of the multivariate Breiman
result in Basrak et al. \cite{basrak:davis:mikosch:2002} shows that
$(X_t)$ is \regvary\ with  index $2\alpha$ and \beao \P(X>x)= 0.5\,
\P(|X|>x) \sim 0.5\, \E|Z|^{2\alpha}\,c_1 \,x^{-2 \alpha}\,. \eeao
Thus both \seq s $(X_t)$ and $(X_t^2)$ are \regvary\ with indices $2
\alpha$ and $\alpha$, respectively. Moreover, {\bf
  (MX)}
is satisfied for both \seq s and we may choose $r_n=n^{\gamma_1}$,
$m_n=n^{\gamma_2}$ for sufficiently small $0<\gamma_1<\gamma_2<1$.
An application of Lemma~\ref{lem:1} yields \eqref{AC}. We omit
details.
\subsubsection{Condition \eqref{TB} for the squared \garch\ process}
Recall the notation $A_t=\alpha_1\,Z_{t-1}^2+\beta_1$,
$B_t=\alpha_0$, $\Pi_t=\prod_{i=1}^t A_i$ and that $(a_n)$ satisfies
$n\,\P(X^2>a_n)\sim 1$. The same arguments as for \eqref{eq:21}
yield \beao X_1^2+\cdots +X_d^2&=& Z_1^2\,\sigma_1^2+\cdots
+Z_d^2\,\sigma_d^2\\[2mm]
&=&\sigma_0^2\,(Z_1^2\,\Pi_1+\cdots +Z_d^2\,\Pi_d)+R_d\,. \eeao
Under the assumption \eqref{eq:24}, $\E (R_d^\alpha)<\infty$, hence
$P(R_d>a_n)=o(P(X^2>a_n))$. This fact and Breiman's result
\cite{breiman:1965} ensure that \beao n\,\P(X_1^2+\cdots +X_d^2>a_n)
&\sim& n\,\P(\sigma_0^2\,(Z_1^2\,\Pi_1+\cdots
+Z_d^2\,\Pi_d)>a_n)\\[2mm]
&\sim & [n\,\E|Z|^{2\alpha}\,\P(\sigma_0^2>a_n)]\,\dfrac{
\E|Z_1^2\,\Pi_1+\cdots
+Z_d^2\,\Pi_d|^\alpha}{\E|Z|^{2\alpha}}\\[2mm]
&\sim &[n\,\P(X^2>a_n)]\,\dfrac{ \E|Z_1^2\,\Pi_1+\cdots
+Z_d^2\,\Pi_d|^\alpha}{\E|Z|^{2\alpha}}\\[2mm]
&\sim &\dfrac{\E|Z_1^2\,\Pi_1+\cdots
+Z_d^2\,\Pi_d|^\alpha}{\E|Z|^{2\alpha}}\\[2mm]
&=&\dfrac{\E|Z_0^2+Z_1^2\Pi_1+ \cdots
+Z_{d-1}^2\,\Pi_{d-1}|^\alpha}{\E|Z|^{2\alpha}}=b_+(d)\,. \eeao In
the last step we used that $A_1$ is independent of $Z_1,\ldots,Z_d$
and that $\E A^\alpha=1$. Write \beao T_d= Z_1^2\Pi_1+ \cdots
+Z_{d}^2\,\Pi_{d}\,. \eeao Observe that \beao T_d\le
\alpha_1^{-1}\left[\Pi_2+\cdots +\Pi_{d+1}\right]\,. \eeao The same
argument as in Section~\ref{subsec:32} proves that the \rhs\
converges a.s. to a finite limit. Hence $T_d\stas T_\infty$ for some
finite limit $T_\infty= \sum_{t=1}^\infty Z_t^2\,\Pi_t\,. $ If
$\alpha\in (0,1]$, we have by concavity of the \fct\
$f(x)=x^\alpha$, $x>0$, \beao
\E[T_{d+1}^\alpha-T_d^\alpha]=\E[(Z_0^2+T_d)^\alpha-T_d^{\alpha}]\le
\E|Z|^{2\alpha}<\infty\,. \eeao If $\alpha\in (1,2)$ we have by the
mean value theorem for some $\xi \in (0, Z_0^2)$ and using the
concavity of the \fct\ $f(x)=x^{\alpha-1}$, $x>0$, \beao
\E[(Z_0^2+T_{d})^\alpha- T_d^\alpha]&=& \alpha \E[( T_d+\xi
)^{\alpha-1}]\\[2mm]
&\le&  \alpha\left[ \E (T_d^{\alpha-1})+\E(|Z|^{2(\alpha-1)})\right]\\[2mm]
&\le & \alpha \,\E(|Z|^{2(\alpha-1)})\left[1+\E(
A^{\alpha-1})+\cdots
  +(\E (A^{\alpha-1}))^{d}\right]\\[2mm]
&=& \alpha \,\E(|Z|^{2(\alpha-1)})\,(1-\E
A^{\alpha-1})^{-1}<\infty\,. \eeao An application of Lebesgue
dominated \con\ yields in the general case $\alpha\in (0,2)$ that
\beao c_+&=&\lim_{d\to\infty}[b_+(d+1)-b_+(d)]=\lim_{d\to\infty}
\dfrac{\E[(Z_0^2+T_{d})^{\alpha}-T_d^{\alpha}]}{\E|Z|^{2\alpha}}\\[2mm]
&=&\dfrac{\E[(Z_0^2+T_\infty)^\alpha-T_\infty^\alpha]}{\E|Z|^{2\alpha}}\in
(0,\infty)\,. \eeao

\subsubsection{Condition \eqref{TB} for the \garch\ process}
Next we calculate the corresponding value $c_+$ for the \garch\
\seq\ $(X_t)$. By the assumed symmetry of $Z$, we have $c_+=c_-$.
Slightly abusing notation, we use the same symbols $b_\pm(d)$,
$(a_n)$, $T_d$, etc., as for $(X_t^2)$. We  choose $(a_n)$ \st\
$n\,\P(|X|>a_n)\sim 1$. We have \beao S_d&=& Z_1\,\sigma_1 +\cdots+
Z_d \,\sigma_d\,. \eeao Since $\E|Z|^{2\alpha}<\infty$ we have for
any $\epsilon>0$, \beao
&&n\,\P(\|(Z_1\sigma_1,\ldots,Z_d\sigma_d)-(Z_1\sigma_1,Z_2|A_2|^{0.5}\sigma_1,\ldots,
Z_d|A_{d}|^{0.5}\sigma_{d-1}\|>\epsilon \,a_n)\\[2mm]
&\le & n\,\P\Big(\sqrt{\alpha_0} \Big(\sum_{i=2}^d
    Z_i^2\Big)^{1/2}>\epsilon\,a_n\Big)\to 0\,.
\eeao Hence (see Jessen and Mikosch \cite{jessen:mikosch:2006},
Lemma 3.12) \beao n\,\P(S_d>a_n)\sim
n\, \P(Z_1\sigma_1+Z_2 A_2^{0.5}\sigma_1+Z_3A_3^{0.5}\sigma_2\cdots+
Z_{d}A_d^{0.5}\sigma_{d-1}>a_n)\,. \eeao Proceeding by induction,
using the same argument as above and in addition Breiman's result,
and writing $\Pi_{s,t}=\prod_{i=s}^t A_i$ for $s\le t$ and
$\Pi_{st}=1$ for $s>t$, we see that \beao
n\,\P(S_d>a_n) &\sim& n\, \P(\sigma_1\,(Z_1+A_2^{0.5}\,Z_2+\cdots + \Pi_{2,d}^{0.5}\,Z_d)>a_n)\\[2mm]
&\sim & [n\, \P(\sigma_1>a_n)]\,\E[(Z_1+A_2^{0.5}\,Z_2+\cdots + \Pi_{2,d}^{0.5}Z_d)_+^{2\alpha}]\\[2mm]
&\sim & \dfrac{\E[(Z_1+A_2^{0.5}\,Z_2+\cdots + \Pi_{2,d}^{0.5}
Z_d)_+^{2\alpha}]}{\E|Z|^{2\alpha}}\\[2mm]
&=& \dfrac{\E[|Z_1+A_2^{0.5}\,Z_2+\cdots + \Pi_{2,d}^{0.5}
Z_d|^{2\alpha}]}{2\,\E|Z|^{2\alpha}} =b_+(d)\,. \eeao In the last
step we used the symmetry of the $Z_t$'s. Writing $
T_d=\sum_{i=1}^{d} Z_i\Pi_{2i}^{0.5}\,, $ we have \beao b_+(d)=
\dfrac{\E [|T_d|^{2\alpha}]}{2\E|Z|^{2\alpha}}\,, \eeao and
\beqq\label{eq:30} |T_d|\le  \sum_{i=1}^\infty
|Z_i|\,\Pi_{2,i}^{0.5}\le \alpha_1^{-1/2}
\sum_{i=1}^\infty\Pi_{2,i+1}^{0.5}\,. \eeqq Since $\E[(\log
A)^{1/2}]=0.5 \E\log A<0$, the \rhs\ converges a.s. to a finite
limit. By a Cauchy \seq\ argument, \beao T_d\stas
T_\infty=\sum_{i=1}^\infty Z_i\,\Pi_{2,i}^{0.5} \eeao for some a.s.
finite $T_\infty$.

Let $(Z_1',A_2')$ be an independent copy of $(Z_1,A_2)$, independent
of  $T_d$. Assume $2\alpha\in (0,1]$. Then by symmetry of the
$Z_t$'s and since $\E[(A_2')^{\alpha}]=1$, using the concavity of
the \fct\ $f(x)=x^{2\alpha}$, $x>0$,
\begin{eqnarray*}
\lefteqn{\E\big[\big| |T_{d+1}|^{2\alpha}-|T_d|^{2\alpha}\big|\big]}\qquad\\
&=& \E\big[\big||Z_1'+(A_2')^{0.5}\,T_d|^{2\alpha}-
|(A_2')^{0.5}T_d|^{2 \alpha}\big|\big]\\[2mm]
&=& \E\big[\big||Z_1'|+(A_2')^{0.5}\,T_d|^{2\alpha}-
|(A_2')^{0.5}T_d|^{2\alpha}\big|\big]\\[2mm]
&=&\E[(|Z_1'|+(A_2')^{0.5}\,(T_d)_+)^{2\alpha}-((A_2')^{0.5}(T_d)_+)^{2\alpha}]+\\[2mm]
&&\E\big[\big|
(|Z_1'|-(A_2')^{0.5}\,(T_d)_-)^{2\alpha}-((A_2')^{0.5}(T_d)_-)^{2\alpha}\big|
I_{\{(A_2')^{0.5}(T_d)_-\le|Z_1'|\}}\big]+\\[2mm]
&&\E\big[\big|((A_2')^{0.5}\,(T_d)_--|Z_1'|)^{2\alpha}-((A_2')^{0.5}(T_d)_-)^{2\alpha}\big|I_{\{(A_2')^{0.5}(T_d)_->|Z_1'|\}}\big]\\[2mm]
&\le &\E|Z|^{2\alpha}\,.
\end{eqnarray*}
For $2\alpha\in (1,2)$ we use the same decomposition as above and
the mean value theorem to obtain \beao
\E\big[\big||T_{d+1}|^{2\alpha}- |T_d|^{2\alpha}\big|\big] &\le &
2\alpha\,\E[(A_2')^{0.5}|T_d|+|Z_1'|]^{2\alpha-1}\\[2mm]
&\le &2\alpha\,\E[(A_2')^{0.5}|T_d|]^{2\alpha-1}+2\alpha
\E|Z|^{2\alpha-1}\,. \eeao The \rhs\ is bounded since
$\E|Z|^{2\alpha}<\infty$ and, using \eqref{eq:30} and
$\E[A^{\alpha-0.5}]<1$, \beao \E|T_d|^{2\alpha-1}&\le&
c\,\sum_{i=1}^\infty \E[\Pi_{2,i+1}^{\alpha-0.5}] =
c\,\sum_{i=1}^\infty (\E [A^{\alpha-0.5}])^i<\infty\,. \eeao Now we
may apply Lebesgue dominated \con\ to conclude that the limit \beao
c_+= \lim_{d\to\infty}[b_+(d+1)-b_+(d)]&=&
\dfrac{\E[|Z_1'+(A_2')^{0.5}T_\infty|^{2\alpha}-|(A_2')^{0.5}T_\infty|^{2\alpha}]}{2\E|Z|^{2\alpha}}
\eeao exists and is finite.

Since $T_\infty$ and $Z_1'$ assume positive and negative values we
have to show that $c_+>0$. First we observe that \beqq\label{eq:31}
c_+= \lim_{d\to\infty} d^{-1} \E[|T_{d}|^{2\alpha}]\,. \eeqq
Applying Khintchine's inequality (see Ledoux and Talagrand
\cite{ledoux:talagrand:1991}) conditionally on $(|Z_t|)$, we obtain
for some constant $c_\alpha>0$, all $d\ge 1$, \beao
 \E[|T_{d}|^{2\alpha}]&\ge& c_\alpha\,
\E\Big(\sum_{i=1}^d Z_i^2\Pi_{2,i}\Big)^{\alpha}\\[2mm]
&=&c_\alpha\, \alpha_1^{-\alpha} \E\Big(\sum_{i=1}^d (\alpha_1
  Z_i^2+\beta_1)\Pi_{2,i}-\beta_1\sum_{i=1}^d
  \Pi_{2,i}\Big)^{\alpha}\\[2mm]
&\ge &c_\alpha\, \alpha_1^{-\alpha}\Big[ \E\Big(\sum_{i=1}^d
  \Pi_{2,i+1}\Big)^{\alpha}-\beta_1^{\alpha}\,\E\Big(\sum_{i=1}^d
  \Pi_{2,i}\Big)^{\alpha}\Big]\\[2mm]
&=&c_\alpha\, \alpha_1^{-\alpha}(1-\beta_1^{\alpha})
\E\Big(\sum_{i=1}^d
  \Pi_{2,i}\Big)^{\alpha}\,.
\eeao Now, Remark~\ref{rem:1} and \eqref{eq:31} imply that $c_+>0$.
\subsubsection{Condition $5$}
We assume $Z$ symmetric. Then $\E X=0$ for $2\alpha\in (1,2)$ and
\eqref{CT} holds for $(X_t)$. For $(X_t^2)$, \eqref{CT} cannot be
satisfied and needs special treatment. If $\alpha \in (1,2)$ all
arguments above remain valid when $X_t^2$ is replaced by $\E
(X_t^2)-\E (X^2)$.

We summarize our results for the \garch\ process $(X_t)$ and its
squares. \bpr Let $(X_t)$ be a strictly stationary \garch\ process
with symmetric iid unit variance noise $(Z_t)$. Assume that $Z$ has
a positive density on $\bbr$  and that \eqref{eq:24} holds for some
positive $\alpha$. Then the  following statements hold.
\begin{enumerate}
\item[\rm (1)]
The \seq s $(X_t)$  and $(X_t^2)$ are \regvary\ with indices
$2\alpha$ and $\alpha$, respectively, and both are strongly mixing
with geometric rate. In particular, \beao \P(X>x)\sim  \frac 12
\E|Z|^{2\alpha}\,c_1\,x^{-2\alpha}\,, \eeao where $c_1$ is defined
in \eqref{eq:34}.
\item[\rm (2)]
Assume $2\alpha\in (0,2)$. Then \beao
(c_1\,E|Z|^{2\alpha}\,n)^{-1/(2\alpha)}S_n\std Z_{2\alpha}\,, \eeao
where $Z_{2\alpha}$ is symmetric $2\alpha$-stable with \chf\
$\psi_{2\alpha}(t) = \exp( -|t|^{2\alpha} \chi_{2\alpha}(t, c_+,
c_+))$, where \beao c_+= \dfrac{\E[|Z_0+ |\alpha_1
Z_0^2+\beta_1|^{0.5} T_\infty|^{2\alpha}-
||\alpha_1Z_0^2+\beta_1|^{0.5}T_\infty|^{2\alpha}]}{2\E|Z|^{2\alpha}}\,,
\eeao and \beao T_\infty=\sum_{t=1}^\infty Z_{t}\,
\prod_{i=1}^{t-1}(\alpha_1\,Z_i^2+\beta_1)^{0.5}\,. \eeao
\item[\rm (3)]
Assume $\alpha\in (0,1)\cup (1,2)$. Then
\[
(c_1\,E|Z|^{2\alpha}n)^{-1/\alpha}\Big(\sum_{t=1}^n X_t^2-
b_n\Big)\std \wt Z_\alpha,\] where
\[ b_n=\left\{\barr{ll}0 &\textrm{\em if\,}\quad \alpha\in (0,1)\\
n\,\E (X^2)& \textrm{\em if\,}\quad \alpha\in (1,2)\,,  \earr\right.
\]
and $\wt Z_\alpha$ is $\alpha$-stable with \chf\
$$\psi_{\alpha}(t) = \exp( -|t|^{\alpha} \chi_{\alpha}(t, c_+, 0)),$$ where
\beao c_+= \dfrac{ \E[(Z_0^2+\wt T_\infty)^\alpha-\wt
T_\infty^\alpha]}{\E|Z|^{2\alpha}}\,. \eeao In the above relation,
\beao \wt T_\infty=\sum_{t=1}^\infty Z_{t+1}^2\,\prod_{i=1}^{t}
(\alpha_1Z_i^2+\beta_1)\,. \eeao
\end{enumerate}
\epr \bre\label{bart} For ARCH(1) processes the above technique of
identification of parameters of the limiting law was developed in
\cite{bartkiewicz:2007}. \ere
\subsection{Stable stationary \seq }
In this section we consider a strictly stationary symmetric
$\alpha$-stable (\sas ) \seq\ $(X_t)$, $\alpha\in (0,2)$, having the
integral \rep\ \beao X_n=\int _E f_n(x)\,M(dx)\,,\quad n\in \bbz \,.
\eeao Here $M$ is an \sas\ random \ms\ with control \ms\ $\mu$ on
the $\sigma$-field ${\mathcal E}$ on $E$ and $(f_n)$ is a suitable
\seq\ of deterministic \fct s $f_n\in L^\alpha(E,{\mathcal E},\mu)$.
We refer to Samorodnitsky and Taqqu \cite{samorodnitsky:taqqu:1994}
for an encyclopedic treatment of stable processes and to Rosi\'nski
\cite{rosinski:1995} for characterizing the classes of stationary
$(X_t)$ in terms of their integral \rep s.

Then for some \sas\ \rv\ $Y_\alpha$, \beqq\label{eq:distr}
S_n=\int_E(f_1(x)+\cdots +f_n(x)) M(dx)\eqd Y_\alpha\,
\Big(\int_E|f_1(x)+\cdots
+f_n(x)|^\alpha\,\mu(dx)\Big)^{1/\alpha}\,. \eeqq Since
$P(Y_\alpha>x)\sim 0.5 c_0 x^{-\alpha}$ for some $c_0>0$ (see Feller
\cite{feller:1971}), we have with $n\,P(|X|>a_n)\sim 1$, \beao
n\,P(S_d>a_n)\sim \dfrac{\int_E|f_1(x)+\cdots
  +f_d(x)|^\alpha\,\mu(dx)}{\int_E |f_1(x)|^\alpha\mu(dx)}=b_+(d)\,,
\quad d\ge 1\,. \eeao Moreover, it follows from \eqref{eq:distr}
that $a_n^{-1} S_n\std Z_\alpha$ for some $Z_\alpha$ \fif\
\beqq\label{eq:last} n^{-1} b_+(n) \to c_+ \eeqq for some constant
$c_+$ and the limit $Z_\alpha$  is \sas , possibly zero.

Since we know the \ds\ of $a_n^{-1} S_n$ for every fixed $n$ we do
not need Theorem~\ref{thm:clt} to determine a \sas\ limit. In the
examples considered above we are not in this fortunate situation. In
the \sas\ case we will investigate which of the conditions in
Theorem~\ref{thm:clt} are satisfied in order to see how restrictive
they are. Since the \fidi s of $(X_t)$ are $\alpha$-stable, {\bf
(RV)} is satisfied. Conditions \eqref{CT} and $\E X=0$ for
$\alpha\in (1,2)$  are automatic. Under \eqref{eq:last}, using the
special form of the \chf\ of a \sas\ \rv , {\bf (MX)} holds for any
\seq\ $m_n\to\infty$. Condition \eqref{AC} is difficult to be
checked. In particular, it does not seem to be known when $(X_t)$ is
strongly mixing. An inspection of the proof of Lemma~\ref{ballou},
using the particular form of the \chf s of the \sas\ \rv s, shows
that \eqref{AC} can be replaced by \eqref{TB} which implies
\eqref{eq:last}. Thus  \eqref{TB} is the only additional restriction
in this case.

\subsubsection*{Acknowledgments}
Parts of this paper were written when Thomas Mikosch visited the
Faculty of Mathematics and Computer Science at Nicolaus Copernicus University Toru\'n
and Universit\'e de Paris-Dauphine and Olivier Wintenberger visited
the Department of Mathematics at the University of Copenhagen. Both
authors take pleasure in thanking their host institutions for
excellent hospitality and financial support. { We would like to thank the
referee for bringing the paper by
Gou\"ezel \cite{gouezel:2004} to our attention. Remark \ref{rem:eps} is also due to the referee.


\end{document}